\documentclass{svmult}
\usepackage{amssymb}
\usepackage[debugshow]{graphicx}
\usepackage{makeidx}
\usepackage{multicol}

\begin{document}
\titlerunning{SPIRAL CHAIN COLORING}
\authorrunning{CAHIT}

\title{A Unified Spiral Chain Coloring Algorithm for Planar Graphs}
\author{I. Cahit}
\institute{Near East University}
\maketitle
\date{}

\begin{abstract}
In this paper we have given a unified graph coloring algorithm for planar
graphs. The problems that have been considered in this context respectively,
are vertex, edge, total and entire colorings of the planar graphs. The main
tool in the coloring algorithm is the use of spiral chain which has been
used in the non-computer proof of the four color theorem in 2004. A more
precies explanation of the proof of the four color theorem by spiral chain
coloring is also given in this paper. Then we continue to spiral-chain
coloring solutions by giving the proof of other famous conjectures of
Vizing's total coloring and planar graph conjectures of maximum vertex
degree six. We have also given the proof of a conjecture of Kronk and
Mitchem that any plane graph of maximum degree $\Delta $ is entirely $%
(\Delta +4)$-colorable. The last part of the paper deals with the three colorability of planar graphs under the spiral chain coloring. We have given an efficient and short proof of the Gr\"{o}tzsch's Theorem that triangle-free planar graphs are $3$-colorable.
\end{abstract}

\section{Introduction}

Without doubt the root of all graph coloring problems e.g., see for example
[13],[33] go to the famous four color map coloring problem of Guthrie [14]
and its solution that is becoming a theorem has a long and strange story
[1]-[3]. Also the lengthy and computer-aided proof(s) and verification of
its correctness by another computer program makes the problem even more
attractive and interesting [1],[17],[34].

The author has given an algorithmic proof to the four color theorem which is
not rely on a computer program but it is based on graph theory notions such
as vertices, edges, cycles, planar graphs etc., in 2004 [19]. The only new
concept introduced in the proof, is a special path in the planar graph
called spiral chain. By using spiral chains and spiral-chain coloring in the
planar graphs proofs have been proposed for several open coloring problems
[20],[21]. The purpose of this paper is to show that spiral chain coloring
algorithm, for at least for planar graphs can be unified to other coloring
problems. In order to show the ability of the spiral-chain coloring
algorithm we have chosen maximal planar graphs for vertex, edge, total and
entire coloring problems. Note that solutions of some of these problems are
still open.

Throughout the paper let us assume that from a plane graph we understand a
maximal planar graph embedded in the plane without crossing of the edges.
The four color problem is to color the vertices of a plane graph with only
four colors so that adjacent vertices receive different colors. The four
color theorem says that four colors is enough for any plane graph. But the
answer of the question of three colorability of the plane graph is open and
only partial results exits, such as planar graphs without triangles or
planar graphs with even triangulations have been shown to be 3-colorable
[22],[23],[25]. We will be re-visited the spiral chain coloring solution of
the Steinberg's three colorability problem in the last section of this paper
for some extra justification [21],[24].

Probably edge-coloring of graphs is almost as old as the four color problem
and comes from its equivalent formulation of Tait [6]. That is 4CT is
equivalent to the coloring the edges of any cubic bridgeless planar graph
with only three colors such that any two incident edges receive different
colors. Again spiral chain edge coloring solution to this problem has been
given by Cahit without relying on the proof of the four color theorem in 2005 [19]. 
The proof is based on the spiral chains of a bridgeless cubic graph and 
coloring them with three colors say Green, Yellow and Red, where Green has priority
over Yellow and Red and Yellow has priority over color Red. Finally in case of color conflict at
two incident edges use appropriate backward Kempe-chain switching to resolve the conflict [19]. 
Tait's coloring is the first example but   
the real starting point of edge coloring of graphs is the famous theorem of
Vizing that states that any graph $G$ has edge-chromatic number $\chi
^{^{\prime }}(G)=\Delta (G)$ or $\Delta (G)+1$, where $\Delta (G)$ is the
maximum vertex degree in $G$ and by "chromatic number" we mean minimum
number of colors. The main problem in edge coloring is to determine which of
these two possibilities holds for a given graph $G$. The graph $G$ is called
it is in Class I if $\chi ^{^{\prime }}(G)=\Delta (G)$ and is in Class II
otherwise. It is a famous conjecture of Vizing that planar graphs with $%
\Delta (G)\geq 6$ are in Class I [16]. Except the case $\Delta (G)=6$ all
other cases have been settled [32].

The problem of simultaneously coloring sets of elements of a graph posed by
Ringel in 1960 [30] who conjectured that the vertices and faces of a plane
graph may be colored with six colors. This has been settled by Borodin in
[29]. Vizing conjectured that the vertices and edges of any graph may be
colored with $\Delta (G)+2$ colors (known as total coloring of graphs) [10].
Similarly a conjecture of Melnikov for edge-face coloring has been settled
by Sanders and Zhao [31].

Lastly plane graph coloring can be considered in its most general forum as
coloring all elements (vertices, edges and faces) simultaneously. This type
coloring has been considered under the name "entire coloring" by Kronk and
Mitchem in 1972 [16]. They also conjectured that any plane graph of maximum
degree $\Delta $ can be colored with $\Delta (G)+4$ colors and showed that
this true for $\Delta (G)=3$. Other results on this conjecture are first by
Borodin for $\Delta (G)\geq 12$ and then $\Delta (G)\geq 7$ and finally
improved to $\Delta (G)\geq 6$ by using discharging and non-existence of an
minimal counter example by Sanders and Zhao [18]. The cases $\Delta (G)\in
\left\{ 4,5\right\} $ remain undecided.

In this paper we have shown that all the above coloring problems of the
plane graphs can be settled algorithmically by the use of spiral chain
coloring technique.

\section{Spiral Chains}

Let $G(V,E)$ be a plane graph with vertex set $V$ and edge set $E$. Assume
that $d(v)\geq 3$ for $v\in V$. The outer-cycle $C_{o0}$ of $G$ is a cycle
for which there is no edges of $G$ remain outer-region of $C_{o}$. Note that
$\left\vert V(C_{o})\right\vert =3$ and all faces are triangles since $G$ is
maximal planar. Let $V(C_{o0})=\left\{ v_{a},v_{b},v_{c}\right\} .$ We
define spiral chain(s) $S_{1}$ of $G$ as a (disjoint) path(s) with a
topological property (we mean the spiral shape of the path) as follows. The
path $P_{1}(v_{a},v_{c})=\left\{ v_{a},v_{b},v_{c}\right\} $ where $%
v_{c}v_{a}\notin P(v_{a},v_{c})$ of $C_{o0}$ is a subpath of spiral chain $%
S_{1},$ that is we select edges of $C_{o0}$ starting from $v_{a}$ in
clockwise direction till $v_{c}$. Then we delete the vertices of $%
P(v_{a},v_{c})$ and obtain the subgraph $G_{1}$ of $G$ which is triangulated
but not necessarily maximal since its outer-cycle $C_{o1}$may have length
greater $3$. Let $V(C_{o1})=\left\{ v_{d,1},v_{d,2},...v_{d,k}\right\} $
where vertices labeled in clockwise direction. Let $v_{d,i}$ be the highest
indexed vertex such that $v_{c}v_{d,i}\in E(G)$. Then we write new extended
spiral subpath of $S_{1}$ as

$P(v_{a},v_{i-1})=$ $P_{1}(v_{a},v_{c})\cup l_{1}\cup P_{2}(v_{d,i},v_{d,i-1}),$

where $P(v_{d,i},v_{d,i-1})=\left\{ v_{d,i},v_{d,i+1},...,v_{d,i-1}\right\} $
and we call $l_{1}=\{v_{c}v_{d,i}\}$ the connecting link-edge of the spiral
sub-paths $P_{1}(v_{a},v_{c})$ and $P_{2}(v_{d,i},v_{d,i-1})$. Similar above
we trace other vertices of the subgraphs $G_{2},G_{3},...,G_{k}$ and
obtain spiral chain of $G$ if $n=|V(G)|=\cup _{i=1}^{k}|V(G_{i})|$ which can be expressed as

$S_{1}=$ $P_{1}(v_{a},v_{c})\cup l_{1}\cup P_{2}(v_{d,i},v_{d,i-1})\cup
l_{2}\cup P_{3}(v_{e,i},v_{e,i-1})\cup l_{3}\cup
P_{4}(v_{f,i},v_{f,i-1})\cup ...$ .

If $\cup _{i=1}^{k}|V(G_{i})|<n$ then this means that there is no link-edge
connecting the last spiral sub-path of $G_{k}$ to the next one in $G_{k+1}$.
In this case we choose the closest vertex $u$ to the last vertex $v$ of $%
S_{1}$such that $vu\notin E(G)$. In general two consecutive spiral chains $%
S_{k}$ and $S_{k+1}$ is separated by an maximal outerplanar subgraph $%
G_{k,k+1}$ such that $vu\notin E(G),v\in S_{k},u\in S_{k+1}$. Start the
spiral chain $S_{2}$ from $u$ as described above. Eventually we obtain
vertex disjoint spiral chains $S_{1},S_{2},...,S_{p}$ when all vertices of $%
G $ have been visited. Hence in general we can write the set of spiral
chains of $G$ as

$\mathcal{S}={\Large \cup }_{i=1}^{p}S_{i}={\Large \cup }%
_{i=1}^{p-1}P_{i}l_{i}\cup P_{p}$

Note that if $\ p=1$ then $S_{1}$ is an Hamiltonian path of $G$ and for $p>1$
some of spiral chains may be a isolated vertex.

For a given graph $G$ the set of $\mathcal{S}$ of spiral chains decompose $G$
into nested vertex disjoint spiral chains. Any vertex can belong exactly one
spiral chain. Since $G$ is maximal planar graph its faces must be triangles
(cycle of length three). A face (triangle) in $G$ can be in three types: $%
\alpha ,\beta $ and $\gamma -$triangles.

\begin{definition}
An triangle in $G$ under the spiral decomposition $\mathcal{S}$ is called $%
\alpha $triangle if all its edges are non-spiral edges, $\beta -$triangle if
exactly two of its edges are non-spiral edges and $\gamma -$triangle if only
one of its edge is an non-spiral edge.
\end{definition}

It is not difficult to see that for any spiral decomposition $\mathcal{S}$
there exits at least one $\gamma $-triangle but we can draw graphs without $%
\alpha $-triangles. The proof of the first statement can be seen that there
is a spiral subpath with a link-edge and a suitable non-spiral edge which
form an maximal outerplanar subgraph of $G.$ But any maximal outerplanar
graph has a $\gamma $-triangle. In Figure 1 we have shown an maximal planar
graph and its spiral chain with all $\beta $-triangles except with one $%
\gamma $-triangle (shown in grey in the figure).

\begin{figure}[htp]
\centering
\includegraphics[width=0.65\textwidth]{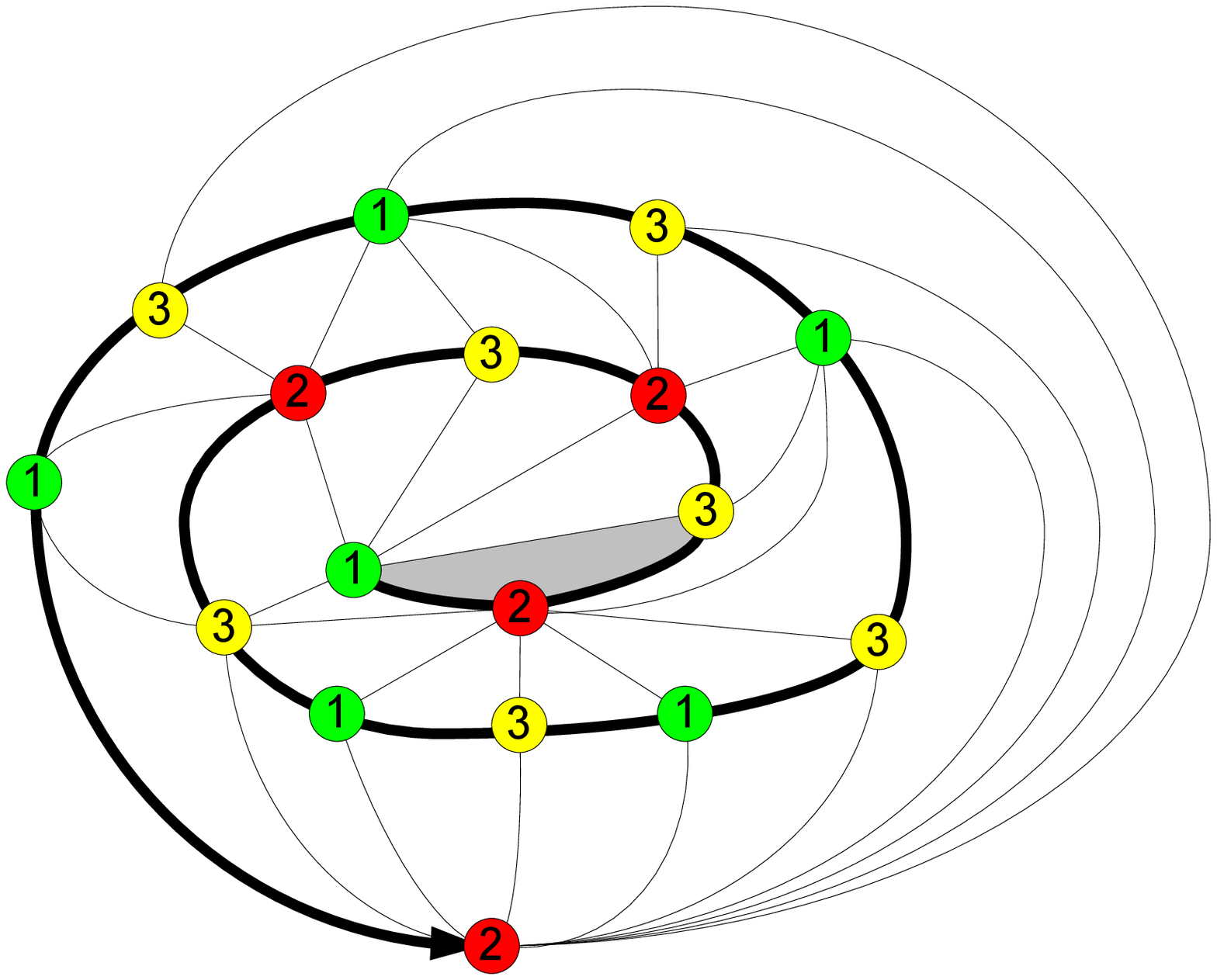}
\caption{Three coloring by spiral chain.}
\label{A:fig:1}       
\end{figure}

\begin{lemma}
Let $\#(\alpha )$ and $\#(\gamma )$ be the number of $\alpha $- and $\ \beta
$-triangles in any spiral chain decomposition $\mathcal{S}$ of $G$. Then $%
\#(\gamma )=\#(\alpha )+1$.
\end{lemma}

Let $G_{o}$ be an maximal outerplanar graph. $G_{o}$ is an triangulated
planar graph in which all vertices are on the (outer-cycle) Hamiltonian
cycle $H_{c}$ and all edges other than the edges of $H_{o}$ can be placed
without crossing into the inside region defined by the $H_{o}.$ Similar
above we can define an triangle in $G_{o}$ as $\alpha $-triangle if all its
edges are non-$H_{o}$ edges, $\beta $-triangle if its two edges are non-$%
H_{o}$ and $\gamma $-triangle if only one of its edge is an non-$H_{o}$ edge.

\begin{lemma}
Let $\#(\alpha )_{o}$ and $\#(\gamma )_{o}$ be the number of $\alpha $- and $%
\ \beta $-triangles in an maximal outerplanar graph $G_{o}$. Then $\#(\gamma
)=\#(\alpha )+2$.
\end{lemma}

Any $(n-1)$ edges of $H_{o}$ form an spiral chain (Hamilton path) in $G$.
Let $e$ be the non-spiral edge of $H_{o}.$ If $e$ was an edge of a $\beta $%
-triangle in $G_{o}$ then becomes an edge of an $\alpha $-triangle in $G$
under $\mathcal{S}$ or similarly if $e$ was an edge of a $\gamma $-triangle
in $G_{o}$ then becomes an edge of a $\beta $-triangle in $G$.

\section{Spiral chain coloring}

In this section we apply spiral chain coloring algorithm to the some of the
planar graph coloring problems which have not completely settled. We have
particularly investigated the undecided cases of the corresponding
conjectures on edge, total and entire coloring conjectures. Let us first
give an complementary justification of the proof of the four color theorem
based on spiral chains of the maximal planar graphs. Assume that $G$ has no
vertex degree smaller than $4$. We will be denoting the colors in several
ways e.g., by numbers 1,2,3,4,.. or by letters $c_{1},c_{2},c_{3},c_{4},...$
or by capital letters $R$ed$,Y$ellow$,G$reen$,B$lue$,....$

\begin{figure}[htp]
\centering
\includegraphics[width=0.65\textwidth]{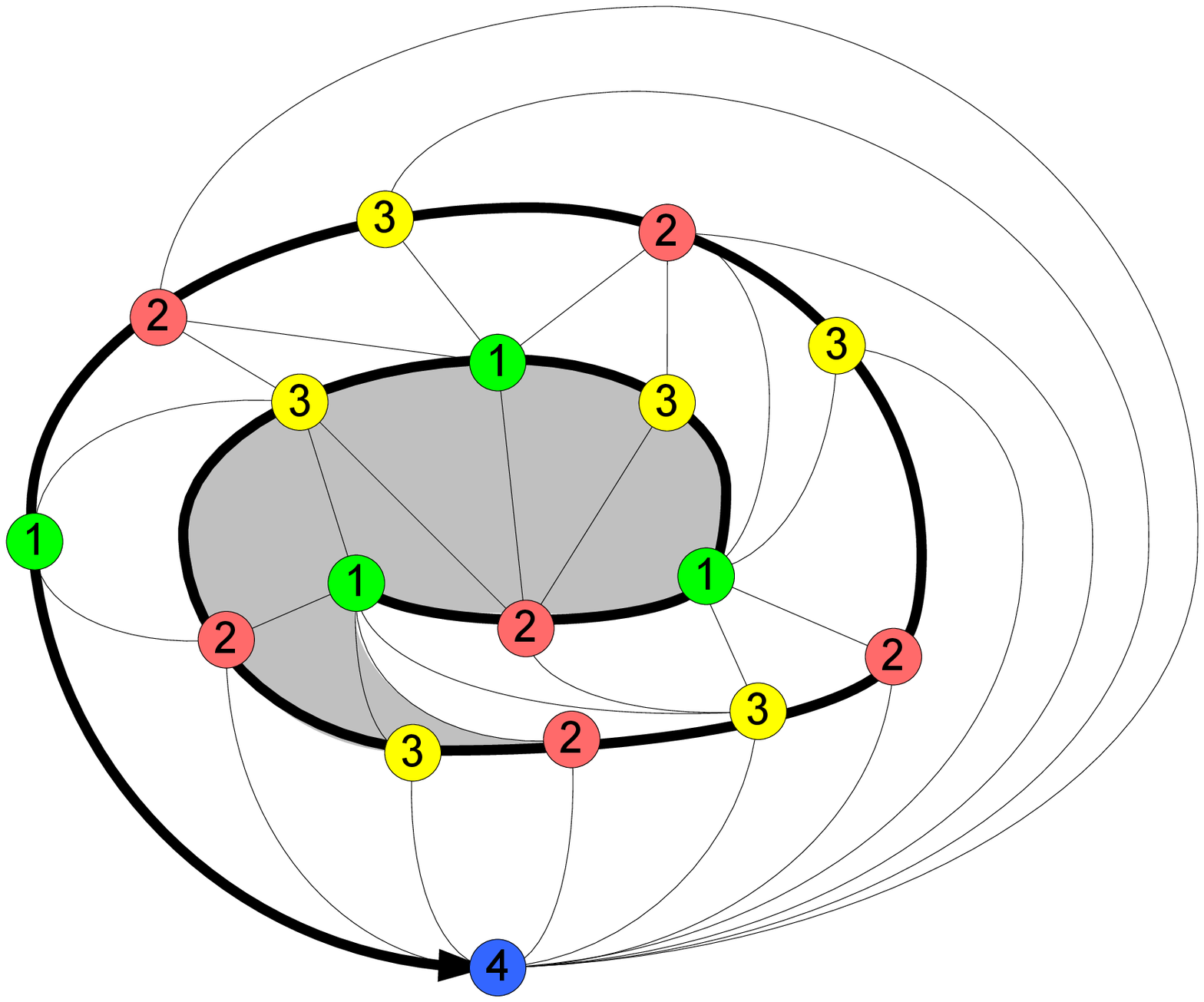}
\caption{\emph{Almost} three colorable graph.}
\label{A:fig:2}       
\end{figure}

\subsection{Spiral Chain Vertex Coloring}

Stockmeyer has shown that 3-colorability is $NP$-complete for planar graphs
but we can decide whether the planar graph is $3$-colorable by using spiral
chain coloring. We call a graph $G$ is \textit{almost} three colorable ($4$%
-chromatic critical) if it is possible to color vertices such that only one
vertex has to be colored by the fourth color. In Fig. 1 we have illustrate $%
3 $-coloring of an maximal planar graph by using spiral chain coloring. The
graph satisfy even triangulation condition of Heawood, so it has to be
3-colorable but we don't need to know that beforehand. In Fig. 2 we have
shown another example of spiral chain coloring which is almost $3$-colorable
in the sense that only once and on the last vertex of the spiral chain the
fourth color is being used.

If $G$ and its spiral chain decomposition has no $\alpha $-triangle then
Figs.1-4 suggest the following lemma:

\begin{lemma}
Any $\alpha $-triangle free maximal planar graph can be $4$-colorable by
spiral chain coloring without Kempe-switch.
\end{lemma}

Let $S_{1},S_{2},...,S_{k}$ be the set of spiral chains of $G.$ Spiral chain
coloring algorithm colors the vertices of $S_{1},S_{2},...,S_{k}$ in reverse
order. Let $V(S_{k})=\{v_{1},v_{2},...,v_{m}\}$ ,$m\leq n$ be the set of
vertices of $S_{k}.$ Spiral segment $S_{k,i}$ of $S_{k}$ is the subset of
vertices $V(S_{k,i})=\{v_{1},v_{2},...,v_{p}\}_{i},i=1,2,...,r,p\leq m\leq n$
such that induced an maximal outerplanar subgraph, where $p$ is the length
of the spiral segment and it is maximum possible with respect to this
property. We also call this the first maximal outerplanar subgraph induced
by the spiral segment vertices as the "\textit{core} of the spiral". That is
spiral segment vertices of the core is $V(S_{k,1})=\{v_{1},v_{2},...,v_{p}%
\}_{1}$. In Fig. 2 the core of the spiral is shown in gray color and the
first spiral segment vertices colored by $G,R,G,Y,G,Y,R,Y,R$ and it can
easily be seen that if we take next vertex (colored by $Y$) to the spiral
segment the core is no longer an maximal outerplanar subgraph. The next
spiral segment together with the previous spiral segment vertices forms
another maximal outerplanar subgraph (in Fig. 2 second spiral segment is $%
Y,R,Y,R,Y,R,G$). Hence we can write the $r$ consecutive spiral segments of $%
S_{k}$ as

$V(S_{k})=V(S_{k,1})\cup V(S_{k,2})\cup ...\cup V(S_{k,r})$ , $%
V(S_{k,1})\cap V(S_{k,2})\cap ...\cap V(S_{k,r})=\phi .$ We also say, for
any three consecutive spiral segments $S_{k,(i-1)},S_{k,i},S_{k,(i-1)}$ of
spiral chain $S_{k}$, spiral segment $\ S_{k(i-1)}$ is a lower-spiral
segment of $S_{k,i}$ and spiral segment $S_{k(i+1)}$ is an upper-spiral
segment of $S_{k,i}$.

In its most general form an maximal planar graph is decomposed into vertex
disjoint spiral chains and each spiral chain is further decomposed, in a
well defined fashion, into vertex disjoint spiral segments. Moreover spiral
chains and spiral segments are in the form of an shelling structure (nested
shells of triangulations). Next we will show that there is an simple and
efficient coloring algorithm that colors the vertices of spiral from the
inner spiral chain toward an outer spiral chain. Since the core spiral
segment $S_{k,1}$ always induce an maximal outerplanar subgraph of the graph
we start coloring core spiral-segment vertices with only three colors, say
green, red and yellow without any color conflict. Let us call to these three
colors as Color Class I, $CC_{1}=\left\{ G,R,Y\right\} $. Note that we have
just three colored the first spiral segment (core) and since it is maximal
outerplanar graph three-coloring with $CC_{1}=\left\{ G,R,Y\right\} $ is
unique.

\begin{itemize}
\item \textit{A }$\beta $\textit{-triangle in the core spiral segment }$%
S_{k,1}.$ That is if start from the very first triangle of the spiral chain
which must be an $\gamma $- or $\beta $-triangle all other triangles induced
by the core vertices uniquely colored by the colors of $CC_{1}$. Now
consider next spiral segment $S_{k,2}$ which upper-spiral segment with
respect to the core segment $S_{k,1}$. We continue the coloring the vertices
of $S_{k,2}$ with the colors in $CC_{1}$ as long as there is no $\alpha $-
or $\gamma $-triangle in the previous core subgraph.

\item \textit{An }$\gamma $\textit{-triangle in the core spiral segment }$%
S_{k,1}$. That is if we have an $\gamma $-triangle with three consecutive
vertices $v_{i-1},v_{i},v_{i+1}\in V(S_{k,1})$ such that $%
(v_{i-1}v_{i}),(v_{i}v_{i+1})\in E(S_{k,1}),$ and $(v_{i-1}v_{i+1})\in E(G)$
then $v_{i-1},v_{i},v_{i+1}$ must use all the distinct colors of $CC_{1},$
say $c(v_{i-1})=G,c(v_{i})=R,c(v_{i+1})=Y$. An vertex $v\in V(S_{k,2})$ that
adjacent to the vertices $v_{i-1},v_{i},v_{i+1}$ forms three $\beta $%
-triangles and force to use the new color $B$ for the vertex $v.$ If $%
(uw)\in $ $S_{k,2}$ such that $(uv_{i}),(wv_{i}),(wv_{i+1})\in E(G)$ are the
non-spiral edges and $c(u)=Y$ then we color $w$ as $c(w)=B$. Hence in these
cases we switch to new three color class $CC_{2}=\left\{ R,Y,B\right\} .$
\end{itemize}

From the above three colors classes $CC_{1}=\left\{ G,R,Y\right\} $ and $%
CC_{2}=\left\{ R,Y,B\right\} $ we define \textit{safe} colors as follows:
The color green $G\in CC_{1}$ is a safe-color with respect to $CC_{2}$ since
$G\notin CC_{2}$ and similarly the color blue $B\in CC_{2}$ is a safe-color
with respect to $CC_{1}$ since $B\notin CC_{1}$. Hence the red $R$ and
yellow $Y$ are non-safe colors for both $CC_{1}$ and $CC_{2}.$ Now let if we
would have no other triangle types (that is no $\alpha $-triangles) in the
spiral decomposition we would color without need of use of Kempe-switch all
spiral segments with the alternating the two three color classes $CC_{1}$
and $CC_{2}.$ That is $S_{k,i}\Longrightarrow CC_{1}=\left\{ G,R,Y\right\}
,i=1,3,5,...$ and $S_{k,i}\Longrightarrow CC_{2}=\left\{ R,Y,B\right\}
,i=2,4,6,...$.

Now consider the last vertex $v_{l}$ of the last spiral segment $S_{k,r}.$
If $r$ the number of spiral segments is odd then color $v_{l}$ as $%
c(v_{l})=G $ and if $r$ is even then color $v_{l}$ as $c(v_{l})=B$. This
completes the proof of the above lemma.

\begin{figure}[htp]
\centering
\includegraphics[width=0.65\textwidth]{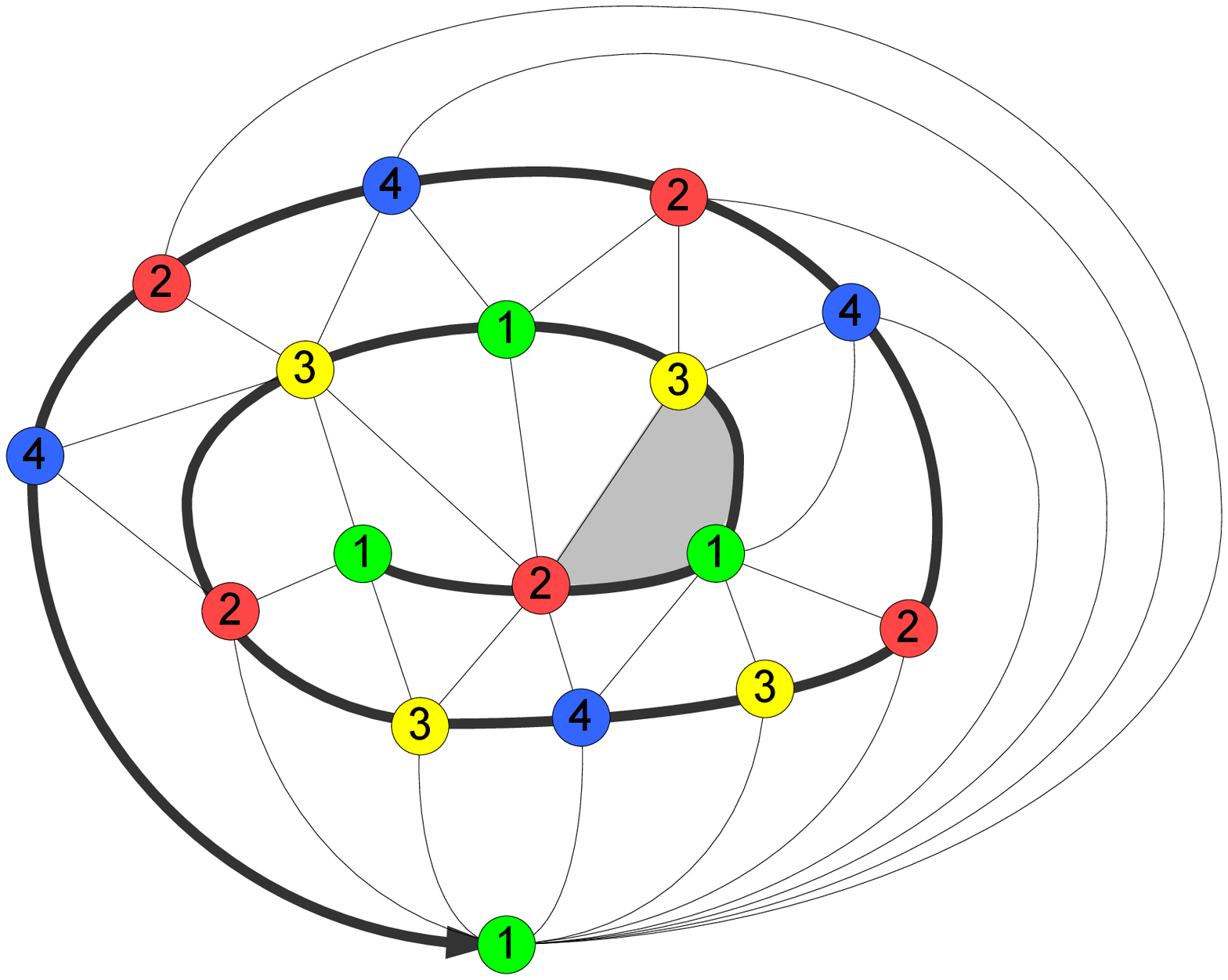}
\caption{Spiral-chain coloring without  Kempe-switch.}
\label{A:fig:3}       
\end{figure}

\begin{figure}[htp]
\centering
\includegraphics[width=0.65\textwidth]{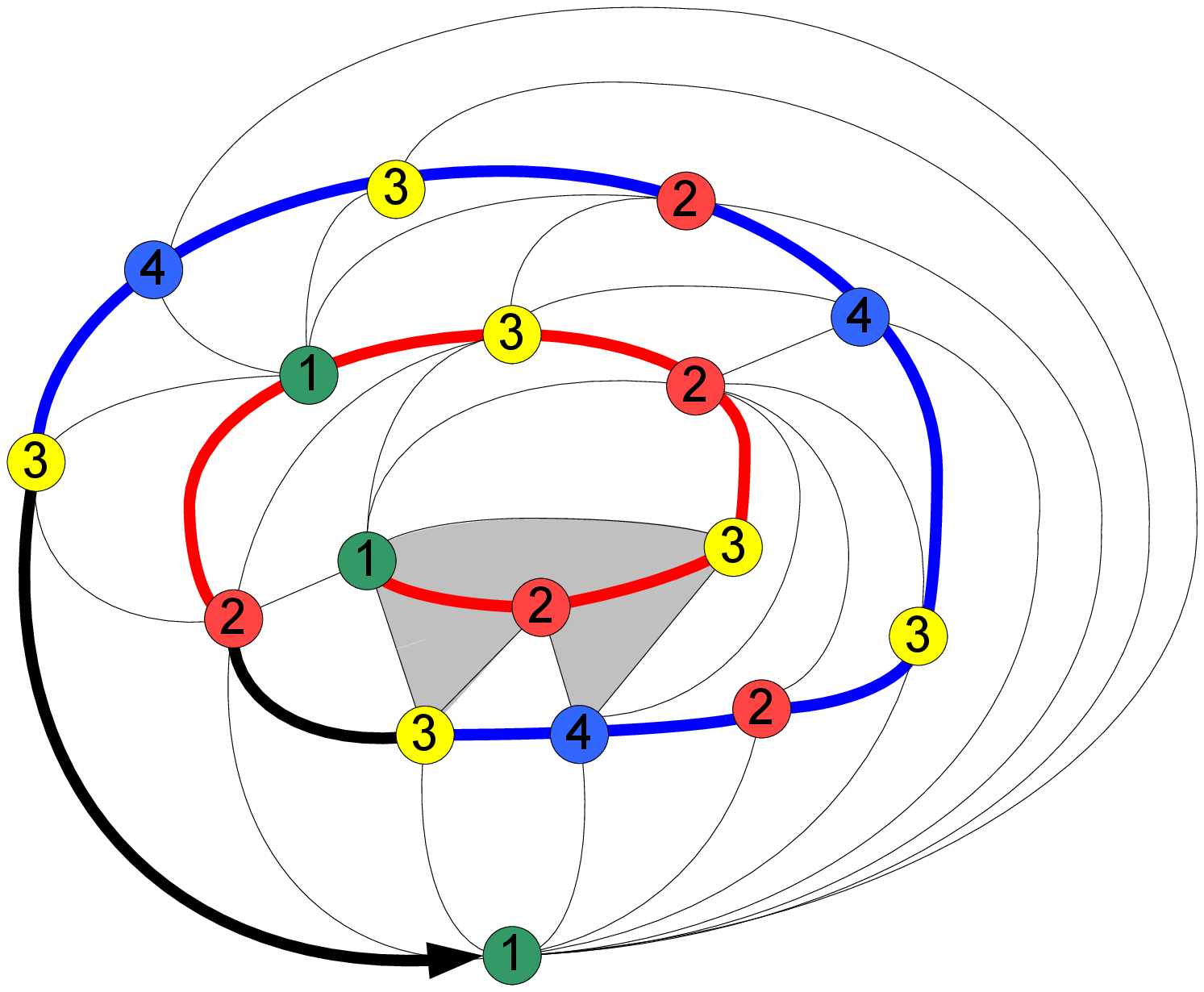}
\caption{Spiral-chain vertex $4$-coloring of an planar graph without an $\alpha$-triangle and with a sailing-boat subgraph.}
\label{A:fig:4}       
\end{figure}

We need a definition of a special subgraph in the spiral chain decomposition
of $G$.

\begin{definition}
Consider two consecutive spiral chains (or spiral segments) $S_{p}$ and $%
S_{p+1}$. Let $v_{i-1},v_{i},v_{i+1}\in S_{p}$ and $v_{r},v_{r+1}\in
S_{p+1}, $ Consider the subgraph formed by $(v_{i-1}v_{i}),(v_{i}v_{i+1})\in
E(S_{p})$ and $(v_{r}v_{r+1})\in E(S_{p+1})$ and the non-spiral edge $%
(v_{i-1}v_{i+1})$ of $S_{p}$ and non-spiral edges $%
(v_{r}v_{i-1}),(v_{r}v_{i}),(v_{r+1}v_{i}),(v_{r+1}v_{i+1}).$ That is the
subgraph formed in this way is called the "sailing boat" and pictorially
looks like a sailing boat in between two consecutive parallel spiral chains
which consist of one $\gamma $-triangle and three $\beta $-triangles. In
other words the sailing boat subgraph is a wheel with five vertices drawn in
the plane like the shape of a sailing boat between two parallel spiral
segments.
\end{definition}

\begin{itemize}
\item \textit{An }$\alpha $\textit{-triangle in the core spiral segment }$%
S_{k,1}.$ It is easy to see that under spiral chains decomposition an
isolate $\alpha $-triangle can create three $\gamma $-triangles. That is in
an triangulation under spiral chain, any sequence of $\beta $-triangles with
a common $\alpha $-triangle edge must end-up with an $\gamma $-triangle.
Recall that all edges of an $\alpha $-triangle are non-spiral edges.
Consider a sailing-boat subgraph between $S_{i-1}$ and $S_{i}$. Let $\left\{
v_{1},v_{2},v_{3},v_{4},....,v_{x},...\right\} \in S_{i-1}$ where $%
v_{1},v_{2},v_{3}$ are the vertices of an $\gamma $-triangle of the sailing
boat and $\left\{ ...,u_{1},u_{2},...\right\} \in S_{i}$. The edge sets of $%
S_{i}$ and $S_{i-1}$ are $(u_{1}u_{2})\in E(S_{i}),\left\{
(v_{1}v_{2}),(v_{2}v_{3}),(v_{3}v_{4})\right\} \in E(S_{i-1})$. The edges
(non-spiral) of the $\alpha $-triangle are $%
(v_{1}v_{3})(v_{3}v_{x})(v_{1}v_{x})$ and $(v_{x}v_{4})\in E(G).$\ Without
loss of generality assume that $S_{i-1}\Longrightarrow CC_{1}=\left\{
G,R,Y\right\} $ and $S_{i}\Longrightarrow CC_{2}=\left\{ B,Y,G\right\} $. \
Spiral chain coloring color the vertices as follows: $%
c(v_{1})=G,c(v_{2})=R,c(v_{3})=Y$ ($\gamma $-triangle in $S_{i-1}$)$%
\Longrightarrow c(v_{x})=R$ ($\alpha $-triangle vertex) and $c(v_{4})=G$ ($%
\beta $-triangle vertex). Now if in $S_{i}$, $c(u_{1})=B$ then we cannot
find proper color for $c(u_{2})=?$ \ But then re-color $c^{\prime }(v_{3})=B$
and $c^{\prime }(u_{2})=Y$ (a single Kempe-switch) to resolve impasse on
vertex $u_{2}$ and maintain $3$-coloring of spiral chain (or segment) $S_{i}%
\footnote{%
Exchange of a safe color of the upper spiral chain with a proper non-safe
color of lower spiral chain can be viewed as preparation the rest of spiral
chain segment vertices for 3-coloring. \ Think of hiding the unwanted
colored spots on the surface of an cake by pushing them with your finger!}$.
Of course spiral chain (or segment) $S_{i-1}$ becomes $4$-coloring since $B$
is a safe color of $CC_{2}$ and $Y$ is a non-safe color.
\end{itemize}

\begin{figure}[htp]
\centering
\includegraphics[width=0.65\textwidth]{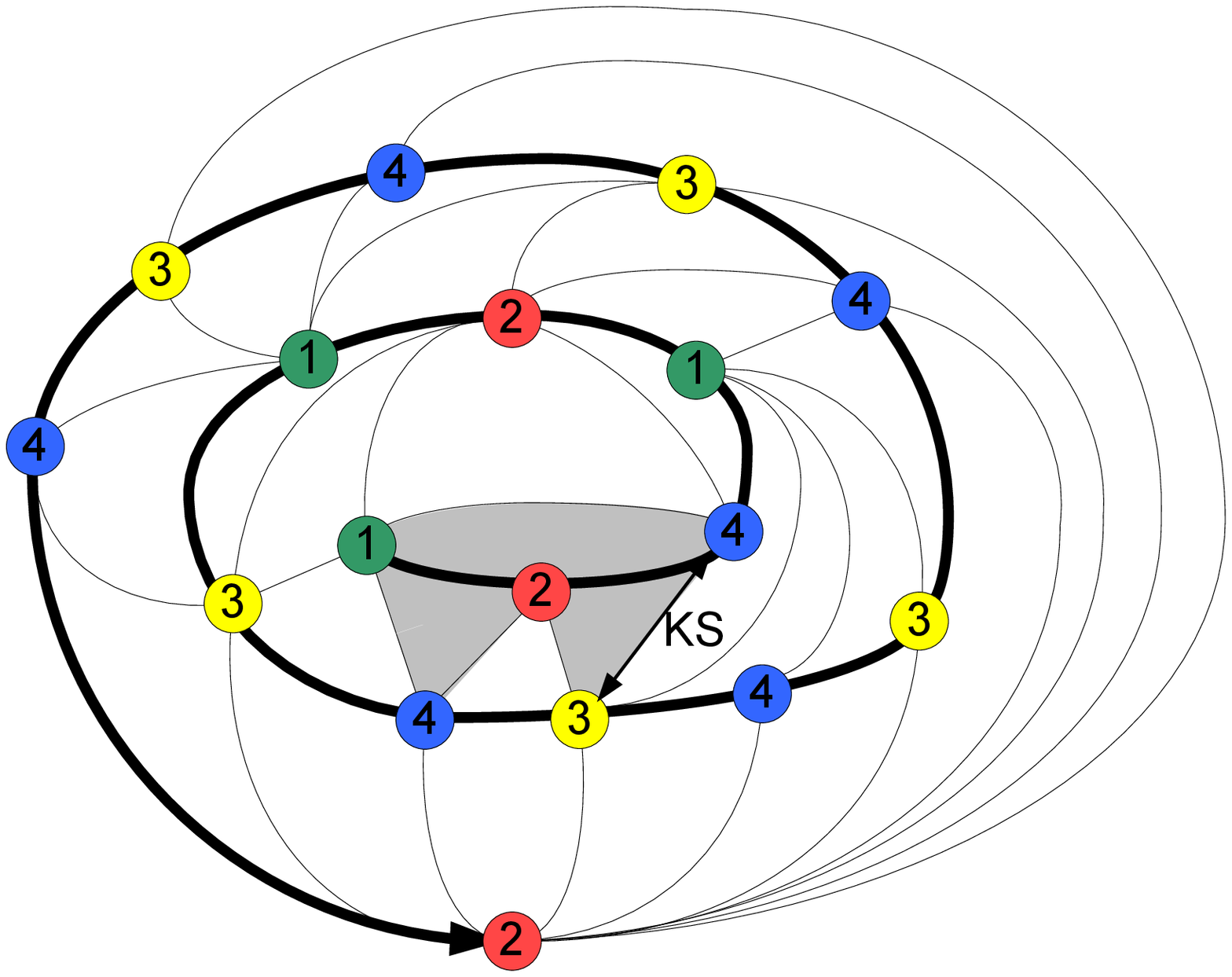}
\caption{Spiral-chain coloring with a Kempe-switch.}
\label{A:fig:5}       
\end{figure}

\textbf{Algorithm 1.}[Description]
Let $S_{1},S_{2},...,S_{k}$ be the set of spiral chains of $G$. Color the
vertices from an inner spiral chain towards an outer spiral chain. Color
spiral chain from inner towards outer spiral segments. Color the core spiral
segment with the color class $CC_{1}=\left\{ G,R,Y\right\} $. For the other
spiral segments use $S_{k,i}\Longrightarrow CC_{1}=\left\{ G,R,Y\right\}
,i=1,3,5,...$ and $S_{k,i}\Longrightarrow CC_{2}=\left\{ R,Y,B\right\}
,i=2,4,6,...$. An vertex in the core-spiral receive an unique color form $%
CC_{1}$ based on the adjacent previously colored triangle. In all spiral
segments other than the core-spiral assign non-safe color to a vertex
whenever is possible. If non-safe color cannot be assigned use respective
safe color of the three color classes. In a spiral segment coloring if an
vertex is in the sailing-boat subgraph and cannot be colored properly then
switch safe color with non-safe color between the parallel spiral segments.
This operation assures three colorability of the current outer spiral
segment at any step. Furthermore three colorability of the outer-spiral
segment assures always to find an safe color to assign to the last vertex of
the spiral chain.

Therefore we write the following theorem from which the four color theorem
follows:

\begin{theorem}
All maximal planar graphs are $4$-colorable.
\end{theorem}

In the next section we will investigate edge-coloring problem of planar
graphs under the spiral chain coloring technique.

\subsection{Spiral Chain Edge Coloring}

In 1979 Seymour has conjectured that there is no planar non-elementary
critical graph [7]. This conjecture implies the four color theorem, the
existence of an algorithm determining the chromatic index of a planar graph
in polynomial time and non-existence of planar class two graph with maximum
degree at least $6$. \ The latter from 1965 Vizing also proved the case $%
\Delta (G)\geq 8$. There are planar class two graphs known with maximum
degree $\ 2,3,4,$ and $5$. The case $\Delta (G)=7$ has been settled by using
discharging method by Gr\"{u}nwald in his Ph.D. thesis in 2000 [8]. The
cyclic spiral coloring algorithm given in this paper not only settles the
case $\Delta (G)=6$ but also answers Seymour's question in affirmative.
Another comment about the cyclic spiral chain coloring algorithm is the
complexity of the $3$-colorability problem of planar graphs.

Let us assume that $G_{a}$ is an \textit{almost} maximal planar graph with
maximum vertex degree $6$ such that all its finite faces are triangles.
Clearly $G_{a}$ may be made fully maximal by joining its all outer-vertices
to another vertex $v_{o}$. An configuration around the vertex $v_{x}\in
S_{i} $ is the subgraph of $G$ induced with all adjacent vertices of $v_{x}$
and the vertex $v_{x}$ itself. Consider three sections of spiral chains $%
S_{i-1},S_{i}$ and $S_{i+1}$. That is spiral-section $S_{i}$ is neighbor
both $S_{i-1}$ and $S_{i+1}$. We say $S_{i+1}$ is upper-spiral neighbor of $%
v_{x}\in S_{i}$ and $S_{i-1}$ is lower-spiral neighbor of $v_{x}\in S_{i}$.
An triangle in between $S_{i}$ and $S_{i+1}$ is called upper triangle $%
\alpha _{u},\beta _{u}$ or $\gamma _{u}$ depending on its type. An triangle
in between $S_{i}$ and $S_{i-1}$ is called lower triangle $\alpha _{l},\beta
_{l}$ or $\gamma _{l}$ depending on its type. Then all triangles of an
configuration with respect to vertex $v_{x}\in S_{i}$ can be written in
anticlockwise direction cyclically as an sequence of triangle types

\begin{center}
$\left\langle \alpha _{l}\beta _{l}\gamma _{l}\alpha _{u}\beta _{u}\gamma
_{u}\right\rangle $
\end{center}

where the order of $\alpha _{i},\beta _{i},\gamma _{i},i=l,u$ depends on the
structure of the configuration but we always start from the first
lower-triangle. For example in Fig.6 triangles around vertex $v_{x}$ are $%
\left\langle \gamma _{l}\alpha _{l}\beta _{l}\beta _{l}\beta _{u}\beta
_{u}\right\rangle $.

\textbf{Algorithm 2.} Spiral Chain Edge Coloring.

Step 1. Find spiral chains $S_{1},S_{2},...,S_{k}$ of $G_{a}.$

Step 2. Color anticlockwise direction the edges incident to the first vertex
of $S_{k}$ starting the spiral-edge. The coloring rule used here is to
assign the first available color from the set of colors $C=\left\{
c_{1},c_{2},...,c_{m}\right\} $.

Step 3. Repeat Step 2 for all edges of $S_{k},S_{k-1},...,S_{1}.$

In Fig.7 we illustrate CSP algorithm for an maximal planar graph with $12$
vertices. The graph has one spiral-chain $S$ and its vertices are $%
v_{12},v_{11},...,v_{1}$. Note that $\ \Delta
(G)=d(v_{9})=d(v_{7})=d(v_{4})=d(v_{1})=6$.

\begin{figure}[htp]
\centering
\includegraphics[width=0.6\textwidth]{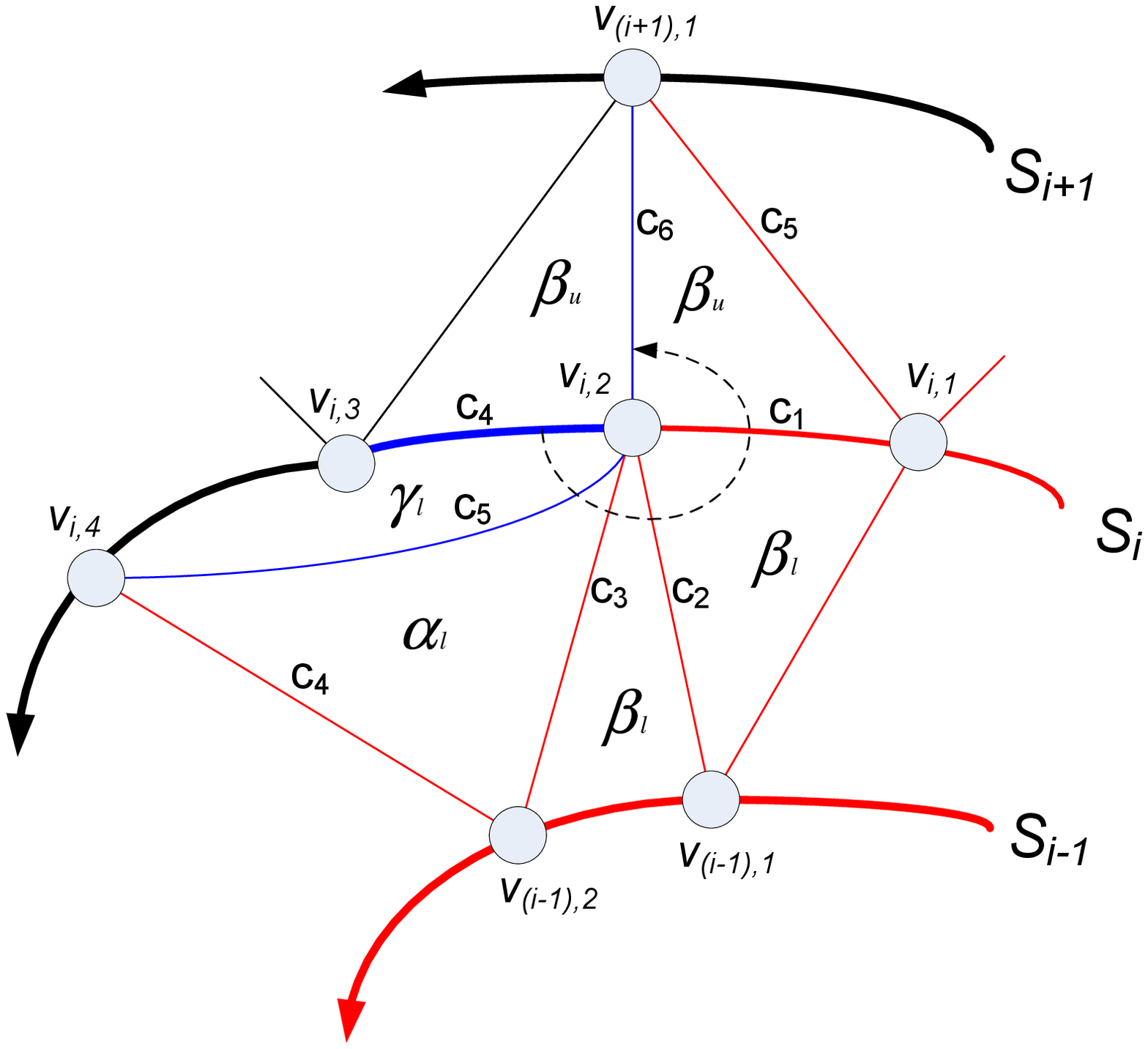}
\caption{Configuration $<\gamma_{l}$ $\alpha_{l}$ $\beta_{l}$ $\beta_{l}$ $\beta_{u}$ $\beta_{u}>$ at vertex $v_{x}$ (red lines already colored, blue lines current and black lines future edges to be colored).}
\label{A:fig:6}       
\end{figure}

\begin{theorem}
The algorithm SCE colors the edges of $G_{a}$ with no more than $\Delta $
colors.
\end{theorem}

\textit{Proof}: Let $C=\left\{ c_{1},c_{2},...\right\} $ be the set of
colors. The proof is based on the argument that cyclic spiral edge coloring
of an configuration at a vertex $v$ never creates an impasse or need of use
more than $\Delta $ colors. Let us list the possible configurations at a
vertex $v$ , where we choose the degree of $v$ as six to show that algorithm
works even in the worst case. Bold lines in counter-clockwise direction
represent three parallel spiral chains $S_{i-1},S_{i},S_{i+1}$.

(1) $\left\langle \beta _{l}\beta _{l}\beta _{l}\beta _{l}\beta
_{l}\right\rangle $

(2) $\left\langle \beta _{l}\beta _{l}\beta _{l}\beta _{l}\beta _{u}\beta
_{u}\right\rangle $

(3) $\left\langle \gamma _{l}\alpha _{l}\alpha _{l}\gamma _{l}\beta
_{u}\beta _{u}\right\rangle $

(4) $\left\langle \beta _{l}\beta _{l}\alpha _{l}\gamma _{l}\beta _{u}\beta
_{u}\right\rangle $

(5) $\left\langle \gamma _{l}\alpha _{l}\beta _{l}\beta _{l}\beta _{u}\beta
_{u}\right\rangle $

(6) $\left\langle \gamma _{l}\beta _{l}\alpha _{l}\beta _{l}\beta _{u}\beta
_{u}\right\rangle $

(7) $\left\langle \beta _{l}\alpha _{l}\beta _{l}\gamma _{l}\beta _{u}\beta
_{u}\right\rangle $

(8) $\left\langle \gamma _{l}\alpha _{l}\beta _{l}\beta _{u}\beta _{u}\beta
_{u}\right\rangle $

(9) $\left\langle \beta _{l}\alpha _{l}\gamma _{l}\beta _{u}\beta _{u}\beta
_{u}\right\rangle $

(10) $\left\langle \gamma _{l}\alpha _{l}\beta _{l}\beta _{u}\alpha
_{u}\beta _{u}\right\rangle $

(11) $\left\langle \beta _{l}\alpha _{l}\gamma _{l}\beta _{u}\alpha
_{u}\beta _{u}\right\rangle $

(12) $\left\langle \beta _{l}\beta _{l}\beta _{l}\beta _{u}\beta _{u}\beta
_{u}\right\rangle $

In configuration $\left\langle \beta _{l}\beta _{l}\beta _{l}\beta _{l}\beta
_{l}\right\rangle $ if $v$ is the last vertex of the spiral chain, since all
edges incident at $v$ have been colored before we terminate the algorithm.
Let us consider configuration (5) $\left\langle \gamma _{l}\alpha _{l}\beta
_{l}\beta _{l}\beta _{u}\beta _{u}\right\rangle $ which is shown in Fig.6.
In Fig.6 edges in "red" \ are already colored, in "blue" are the current
edges and in "black" lines the future edges that to be colored. Let the
three parallel spiral chains be $S_{i-1}=\left\{
...,v_{i-1,1},v_{i-1,2},...\right\} ,S_{i}=\left\{
...,v_{i,1},v_{i,2},v_{i,3},v_{i,4},...\right\} ,S_{i+1}=\left\{
...,v_{i+1,1},...\right\} .$

W.l.o.g. let $%
c(v_{i,2}v_{i,1})=c_{1},c(v_{i,2}v_{i-1,1})=c_{2},c(v_{i,2}v_{i-1,2})=c_{3}.$

If $c(v_{i,4}v_{i-1,2})=c_{4}$ and $c(v_{i+1,1}v_{i,1})=c_{5}$ then put $%
c(v_{i,2}v_{i,3})=c_{4},c(v_{i,4}v_{i,2})=c_{5},$ and $%
c(v_{i,2}v_{i+1,1})=c_{6}.$

If $c(v_{i,4}v_{i-1,2})=c_{5}$ and $c(v_{i+1,1}v_{i,1})=c_{6}$ then put $%
c(v_{i,2}v_{i,3})=c_{4},c(v_{i,4}v_{i,2})=c_{6},$ and $%
c(v_{i,2}v_{i+1,1})=c_{5}.$

If $c(v_{i,4}v_{i-1,2})=c(v_{i+1,1}v_{i,1})=c_{5}$ then put $%
c(v_{i,2}v_{i,3})=c_{5},c(v_{i,4}v_{i,2})=c_{4},$ and $%
c(v_{i,2}v_{i+1,1})=c_{6}.$

If $c(v_{i,4}v_{i-1,2})=c(v_{i+1,1}v_{i,1})=c_{6}$ then put $%
c(v_{i,2}v_{i,3})=c_{6},c(v_{i,4}v_{i,2})=c_{4},$ and $%
c(v_{i,2}v_{i+1,1})=c_{5}.$

Therefore in all cases (coloring of the edges denoted in blue in Fig.6) it
is possible to complete coloring the edges incident $v_{i,2}$ without
needing the seventh color $c_{7}$. It can be verified that this is true for
all other configurations (1)-(12). That is chromatic index $\chi ^{\prime
}(G_{a})=6$ when $\Delta =6$. This completes the proof of the theorem.

\begin{theorem}
The edges of an maximal planar graph $G$ can be colored with $\Delta =6$
colors.
\end{theorem}

\begin{figure}[htp]
\centering
\includegraphics[width=0.7\textwidth]{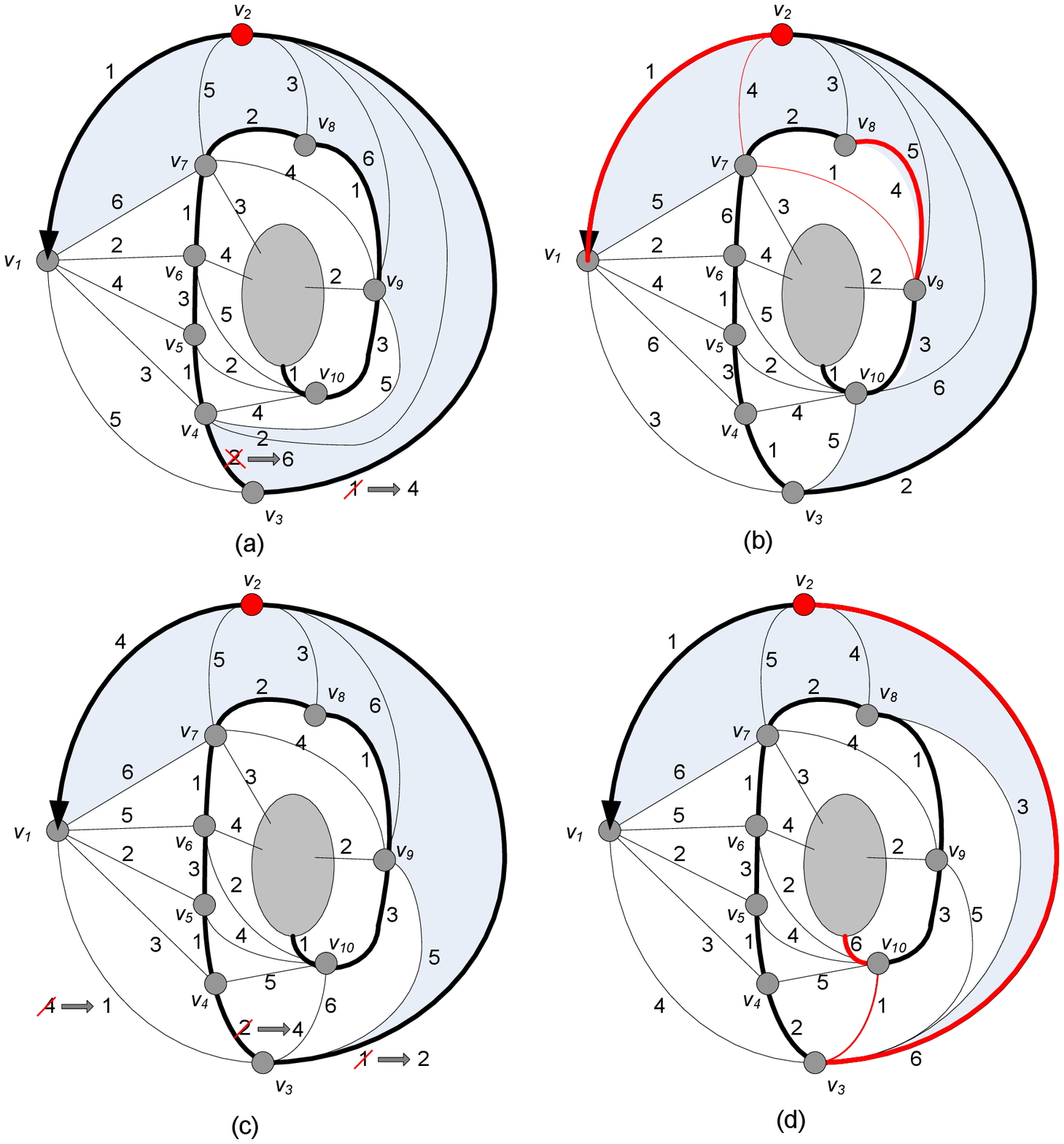}
\caption{Termination of spiral chain edge coloring algorithm with respect to the configuration at vertex $v_2$.}
\label{A:fig:7}       
\end{figure}

\textit{Proof }: We note that outer-cycle of an maximal planar graph has a
length three and $G_{a}$ is a subgraph of $\ G$. Since we assume that $%
\Delta (G)=6$ we can obtain $G_{a}$ by deleting some of the outer vertices
of $G$. Let $S_{1},S_{2},...,S_{k}$ be the spiral chains of $\ G$. By
Theorem 5 we know that spiral chain edge coloring algorithm colors the edges
of $G_{a}$ without any impasse, where in the edge-coloring the order spiral
chains is $S_{k},S_{k-1},...$. When we complete the coloring of $G_{a}$ for
some spiral chain $S_{p}$ we can continue in the same way for the other
edges of spiral chains in $G$ as long as edges are not incident to the last
vertex of the spiral chain. Let $v_{1}v_{2}$ be the last edge of the spiral
chain $S_{1}$. In order to complete the proof we have to give termination
condition of the algorithm for edge coloring of $G$ with no more than $%
\Delta (G)$ colors. Since $G$ is an maximal planar graph with $\Delta (G)=6,$
in the worst case (from the point of the spiral edge coloring) we may assume
that the degree of the first vertex $v_{1}$ of $S_{1}$ as $d(v_{1})=6$. We
also assume that the vertices of the spiral-chains are ordered from an
outer-vertex towards inner vertices of $G$ as $V=\{v_{1},v_{2},...,v_{n}\}.$
That is spiral-chains can be written as $V(G)=S_{1}(V_{1})\cup
S_{2}(V_{2})\cup ...\cup S_{k}(V_{k}),$where $S_{i}(V_{i})\neq \phi $ and $%
S_{1}(V_{1})\cap S_{2}(V_{2})\cap ...\cap S_{k}(V_{k})=\phi $. Clearly $%
v_{1}v_{2}$ is a spiral-chain edge and $v_{1}v_{i},i=3,4,...,7$ are
non-spiral edges. Based on the configuration at vertex $v_{2}$ we consider
the following:

\textit{Case 1.} Configuration $\ \left\langle \beta _{l}\beta _{l}\beta
_{l}\alpha _{l}\gamma _{l}\right\rangle $ at the vertex $v_{2}\ ,$ with $%
\deg (v_{2})=6$. Clearly since $v_{2}$ is an outer-vertex all triangles of
configurations at $v_{2}$ must be lower triangles. The graph $G$ with $%
\left\langle \beta _{l}\beta _{l}\beta _{l}\alpha _{l}\gamma
_{l}\right\rangle $ is shown in Fig.7(a) together with an proper
edge-coloring with $\chi ^{\prime }(G)=6.$ Spiral chain coloring works
without any impasse up to the vertex $v_{3}$. Cyclic coloring of the edges
incident to vertex $v_{4}$ cannot be possible for the edge as $%
c(v_{4}v_{2})=6$ since $v_{9}v_{2}$ has already been colored as $%
c(v_{9}v_{2})=6$ before. Cyclic one step shift of the colors $2,3,6$
respectively on the edges $v_{4}v_{3},v_{4}v_{2},v_{4}v_{1}$ resolves this
impasse. That is $c(v_{4}v_{3})=2,c(v_{4}v_{2})=3,c(v_{4}v_{1})=6$. \ Then
we can put $c(v_{3}v_{2})=4$ and $c(v_{2}v_{1})=1$ and complete spiral edge
coloring of $G.$

\textit{Case 2.} Configuration $\ \left\langle \beta _{l}\beta _{l}\beta
_{l}\beta _{l}\beta _{l}\right\rangle $ at the vertex $v_{2},$with $\deg
(v_{2})=6$.The graph $G$ with $\left\langle \beta _{l}\beta _{l}\beta
_{l}\beta _{l}\beta _{l}\right\rangle $ is shown in Fig. 7(b) together with
an proper edge-coloring with $\chi ^{\prime }(G)=6.$ In this configuration
spiral chain edge coloring algorithm faces with the impasse at the last edge
$v_{2}v_{1}$. That is we cannot assign $c(v_{2}v_{1})\neq 1$ since $%
c(v_{2}v_{7})=1$. But we then have $(4,1)-$Kempe chain $%
(v_{2},v_{7},v_{9},v_{8})$. So we can re-color edges of this Kempe-chain as $%
c(v_{2}v_{7})=4,c(v_{7}v_{9})=1,c(v_{9}v_{8})=4$ and open room for the edge $%
v_{2}v_{1}$ to be colored as $c(v_{2}v_{1})=1$.

\textit{Case 3.} Configuration $\ \left\langle \beta _{l}\beta _{l}\beta
_{l}\beta _{l}\right\rangle $ at the vertex $v_{2},$ with $\deg (v_{2})=5$.
Similar to the \textit{Case 1} above (see Fig.7(c)).

\textit{Case 4}. Configuration $\ \left\langle \beta _{l}\beta _{l}\beta
_{l}\right\rangle $ at the vertex $v_{2},$ with $\deg (v_{2})=4$. Similar to
the \textit{Case 2} above (see Fig.7(d)).

Therefore combining the result that $\chi \prime (G)=\Delta (G)$ if $G$ is
simple, planar and $\Delta (G)\geq 7$ (Vizing 1965 and Sanders and Zhao
2001) and Theorem 6 above we can write:

\begin{theorem}
Planar graph $G$ with $\Delta (G)\geq 6$ is a Class 1.
\end{theorem}

Fig. 10 illustrates spiral edge coloring algorithm for an maximal planar
graph $\ G$ with $\ 12$ vertices. Increasing vertex numbers indicate the
spiral chain. Therefore edge coloring starts from the first edge $%
(v_{12}v_{11})$ and assign colors $C=\{c_{1},c_{2},c_{3},c_{4},c_{5},c_{6}\}$
as cyclically in the counterclockwise direction and continue in this way for
the other spiral chain edges $(v_{11}v_{10}),(v_{10}v_{9}),....,(v_{2}v_{1})$%
. In the figure we have shown colors as red, yellow, green, light-blue, blue
and pink where color red has the highest and color pink has the lowest
priority in the cyclic assignment of the colors. Since degree of vertex $%
\deg (v_{2})=4<\Delta (G)$, configuration at $v_{2}$ would not create any
impasse that would otherwise require Kempe-chain switching in order to
resolve the impasse. Hence spiral edge-coloring confirms that $\chi ^{\prime
}(G)=\Delta $ for $\Delta (G)=6$.

\section{Total and Entire Coloring}

It is easy to reach the conclusion that chromatic number in the vertex
coloring of a graph is related with the maximum size of complete graph minor
e.g., Hadwiger Conjecture and chromatic index in the edge coloring of a
graph is related with the maximum vertex degree. In the total coloring in
which vertices and edges of a graph simultaneously colored and in the entire
coloring of a planar graph in which vertex, edge and faces simultaneously
colored, chromatic numbers are related mainly with the maximum vertex
degree. But as the total and entire colorings of $K_{4}$ in Fig. 8 shows
vertex colors are dominant over edge and face colors in the process of
finding exact colorings. Let us use the terminology used in [22]; denote the
vertex, edge, and face sets of $G$ by $V(G),E(G),$ and $F(G),$ respectively.
An total coloring of $G$ (where $G$ \ may not be a planar graph) is a
function assigning values (colors) to the elements of $V(G)\cup E(G)$ and an
entire coloring of $G$ (here $G$ is necessarily an (plane) planar graph) is
a function assigning values to elements $V(G)\cup E(G)\cup F(G)$ in such a
way that any two distinct adjacent/incident elements receive distinct
colors. The total chromatic number $\chi(G)$ a graph $G$ is the least number of colors needed in any total coloring
of $G$. Total coloring conjecture (Behzad, Vizing) asserts the total
chromatic number of any graph is bounded by $\chi(G)\leq \Delta (G)+2.$ Here we will be dealing with the open case of planar
graphs with maximum vertex degree $\Delta (G)=6$.

\ Let $C=\left\{ c_{1},c_{2},...,c_{k}\right\} =\left\{ 1,2,3,....,k\right\}
$ be the set of colors. We will assume that in assigning an color to the
element of a graph color $c_{i}$ has a priority over color $c_{j}$ if $i<j$.

Denote by the sets $C_{v},C_{e}$ and $C_{f}$ distinct colors used in the
coloring vertices, edges and faces \ of a graph $G$. For example for the
colorings of $K_{4}$ in Fig. 8(a) shows a total coloring of $K_{4}$ with $C_{v}^{t}=\left\{
1,2,3,4\right\} ,C_{e}^{t}=\left\{ 1,2,3,4,5\right\} $ and

Fig. 8(b) shows an entire coloring of $K_{4}$ with $C_{v}^{e}=\left\{
1,5,6,7\right\} ,C_{e}^{e}=\left\{ 1,2,3\right\} ,C_{f}^{e}=\left\{
4,5,6,7\right\} $ where upper subscript denote type of the coloring and
lower subscript denotes type of the element in the graph $G$.

Note that we have $C_{v}^{e}\cap C_{e}^{e}=$\ $\phi ,C_{f}^{e}\cap
C_{e}^{e}=\phi $ hence we may say that for entire coloring of $K_{4}$ edge
colors have no influence on the vertex and face colors.

\begin{figure}[htp]
\centering
\includegraphics[width=0.6\textwidth]{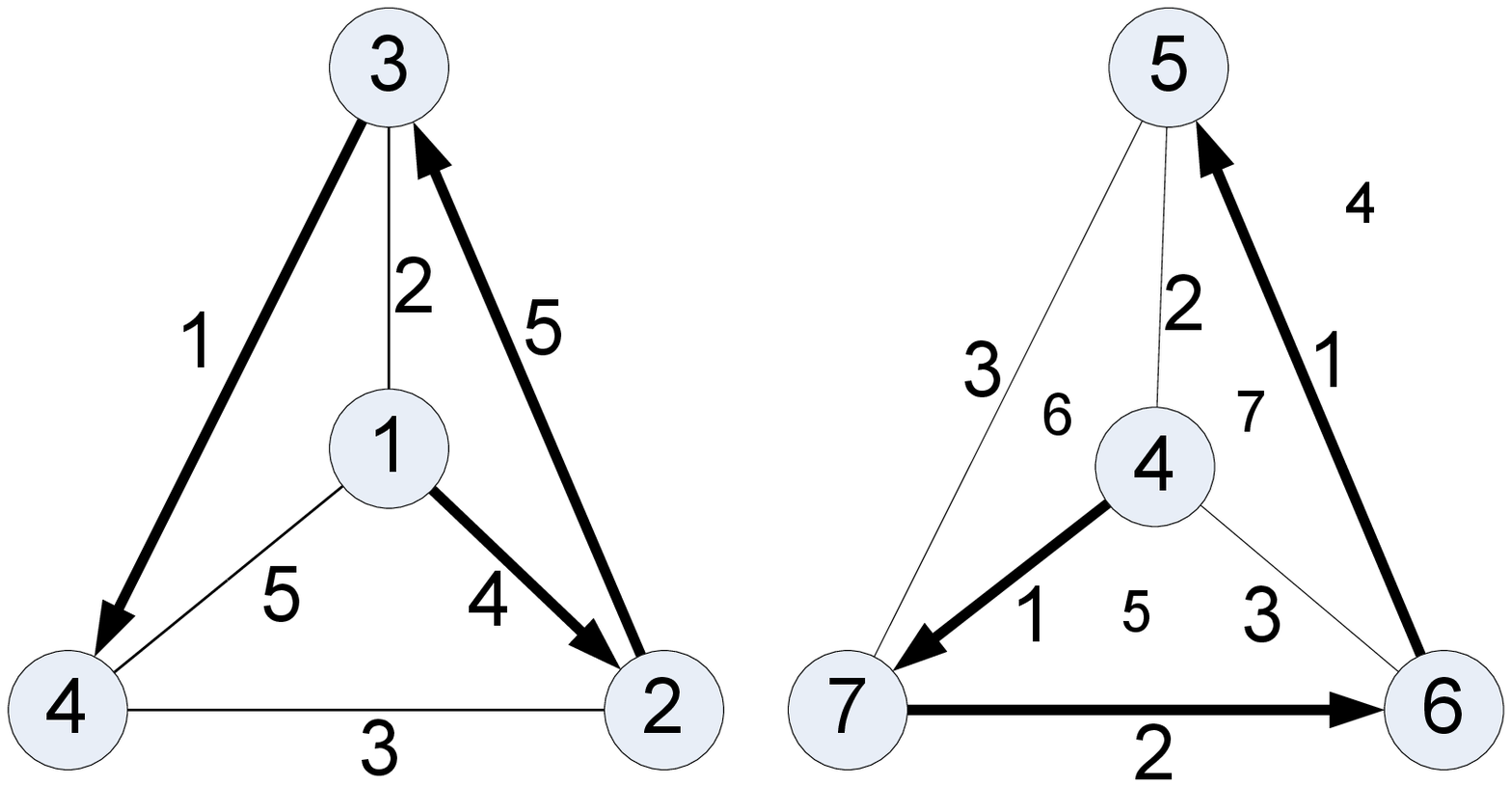}
\caption{Total and entire coloring of $K_4$.}
\label{A:fig:8}       
\end{figure}

\subsection{Spiral Chain Total Coloring}

Kempe chain in an vertex (edge) coloring of a graph is in most general form
is connected two colored subgraph. When the vertices and edges are colored
as in total coloring we can talk about mixed Kempe chain.

\textbf{Definition}
In an total coloring of $G$ for two vertices $v$ and $u$ $(u\neq v)$ if
the edges of an path $P(v,u)$ ,between $v$ and $u$ colored by two colors $
c_{i}$ and $c_{j}$ and $c(u)=c_{j}$ we say the $P(v,u)$ is a (mixed) $m$
-Kempe chain if we have forall $w: (uw)\in E$ and $c(w)\neq c_{i}$.

\textbf{Algorithm 3.}Spiral Total Coloring.
Let $C=\left\{c_{v1},c_{v2},c_{v3},c_{v4},c_{e5},c_{e6},....,c_{e(\Delta +2)}\right\} $ be
the set of colors. Initially we intentionally reserve the first four colors
for the vertices and the other $(\Delta -2)$ colors for the edges of $G$.

Step 1. Color the vertices of $\dot{G}$ by using spiral chain vertex
coloring algorithm with the colors $c_{v1},c_{v2},c_{v3},c_{v4}$.

Step 2. In this step we color the edges of $G$ using spiral chain edge
coloring algorithm using the colors in the set $C$. While assigning an color
to an edge give always priority to low index color.

Step 3. If all edges of $G$ colored with no more than $\Delta +2$ colors
terminate edge-coloring algorithm. If the last edge of the spiral chain
creates color conflict then use $m$-Kempe switch to resolve the color
conflict. This will be explained in detail below.

Let $(v_{1},v_{2},...,v_{n})$ be the set of the vertices of a spiral chain $%
S_{1}$ of $G$ with $\Delta (G)=6$. If $G~$\ has more than one spiral chains
our argument valid for each spiral chains. It is easy to see that at Step 2
of the algorithm as long as we color the edges incident to vertex $v_{i},i>2$
i.e, internal vertex of the spiral chain$,$ since we have colored four edges
form the subset of colors reserved for the edges $C_{e}=\left\{
c_{e5},c_{e6},....,c_{e(\Delta +2)}\right\} $ of $C$ and for the other two
edges incident $v_{i}$ from the subset of colors $C_{v}=\left\{
c_{v1},c_{v2},c_{v3},c_{v4}\right\} $ which has been used originally for the
vertices of $G$. But when we arrive to color the last spiral-chain edge $%
(v_{2}v_{1})$ if $\deg (v_{2})=\Delta (G)$ we may have a situation that $%
c(v_{2}v_{1})=c(v_{2}v_{y})$ where edges incident $v_{y}$ has been handled
before by the algorithm. That is proper coloring of the edges incident to
the last vertex $v_{1}$ of the spiral chain leads to color-conflict. In this
case there must another vertex $v_{x}$ ($\deg (v_{x})<6$) adjacent to $v_{1}$
(that is $(v_{1}v_{x})$ is an non-spiral edge) such that $%
c(v_{1}v_{x})=c_{ei},i\in C_{e}$ and $c(v_{x})=c(v_{2}v_{y}).$ Hence the
four vertices $v_{y},v_{2},v_{1},v_{x}$ form a $m$-Kempe chain. By
performing $m$-Kempe-chain switching we can recolor the edges and the vertex
$v_{x}$ as:

$c^{\prime }(v_{2}v_{y})=c(v_{2}v_{y})=c_{vi}$

$c^{\prime }(v_{2}v_{1})=c(v_{1}v_{x})=c_{e_{i}}$

$c^{\prime }(v_{1}v_{x})=c(v_{x})$

$c^{\prime }(v_{x})=c(v_{1}v_{x})=c_{ei}.$

resolve the edge color conflict and complete total coloring of $\ G$ with $%
\Delta +2$ colors.

From the spiral chain edge and total coloring algorithms we have contributed
to the famous total coloring conjecture of Vizing and Behzad:

\begin{theorem}
The total chromatic number $\chi ^{\prime }(G)$ of planar graphs is $\chi
^{\prime }(G)\leq \Delta (G)+2$.
\end{theorem}

For an illustration consider the total coloring of the graph in Fig.10.
Algorithm completes its vertex and edge coloring by assigning "green" to the
first edge $(v_{2}v_{1})~$of the spiral chain. But we have $%
c(v_{2}v_{7})=c(v_{2}v_{1})$. However we have a $m~$-Kempe chain $%
v_{2},v_{1},v_{3}$. So we apply $m$-Kempe chain switch as follows:

$c(v_{2}v_{1})="pink"\Longrightarrow c^{\ast }(v_{2}v_{1})="green"$

$c(v_{1}v_{3})="green"\Longrightarrow c^{\ast }(v_{1}v_{3})="pink"$

$c(v_{3})="pink"\Longrightarrow c^{\ast }(v_{3})="green".$

So edge color conflict resolved i.e., $c(v_{2}v_{7})\neq c^{\ast
}(v_{2}v_{1})$.

\begin{figure}[htp]
\centering
\includegraphics[width=0.6\textwidth]{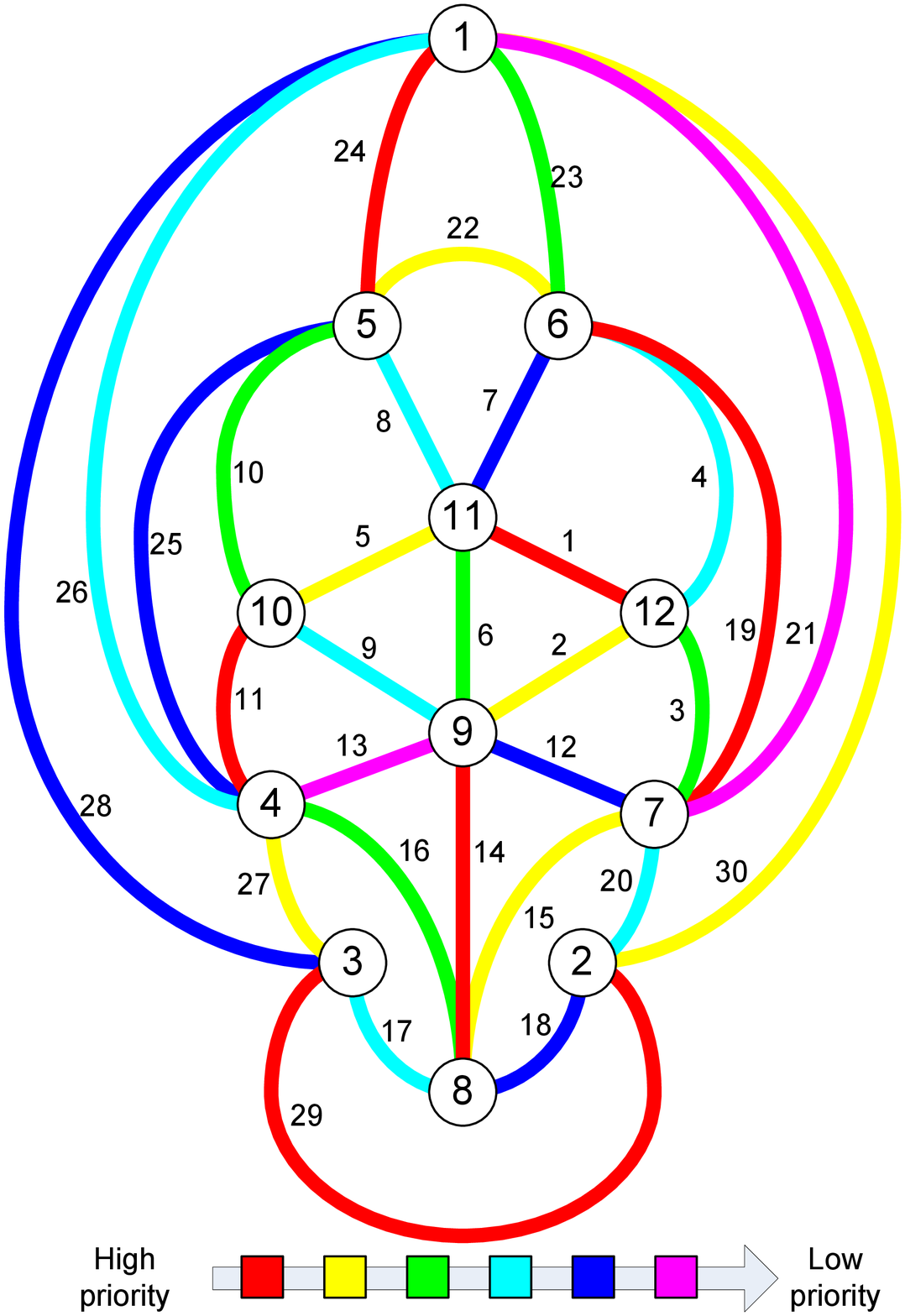}
\caption{Spiral chain edge-coloring of a graph with $\Delta=6$.Increasing vertex numbers indicate the spiral chain. Increasing edge numbers indicate order of edge coloring.}
\label{A:fig:9}       
\end{figure}

\begin{figure}[htp]
\centering
\includegraphics[width=0.6\textwidth]{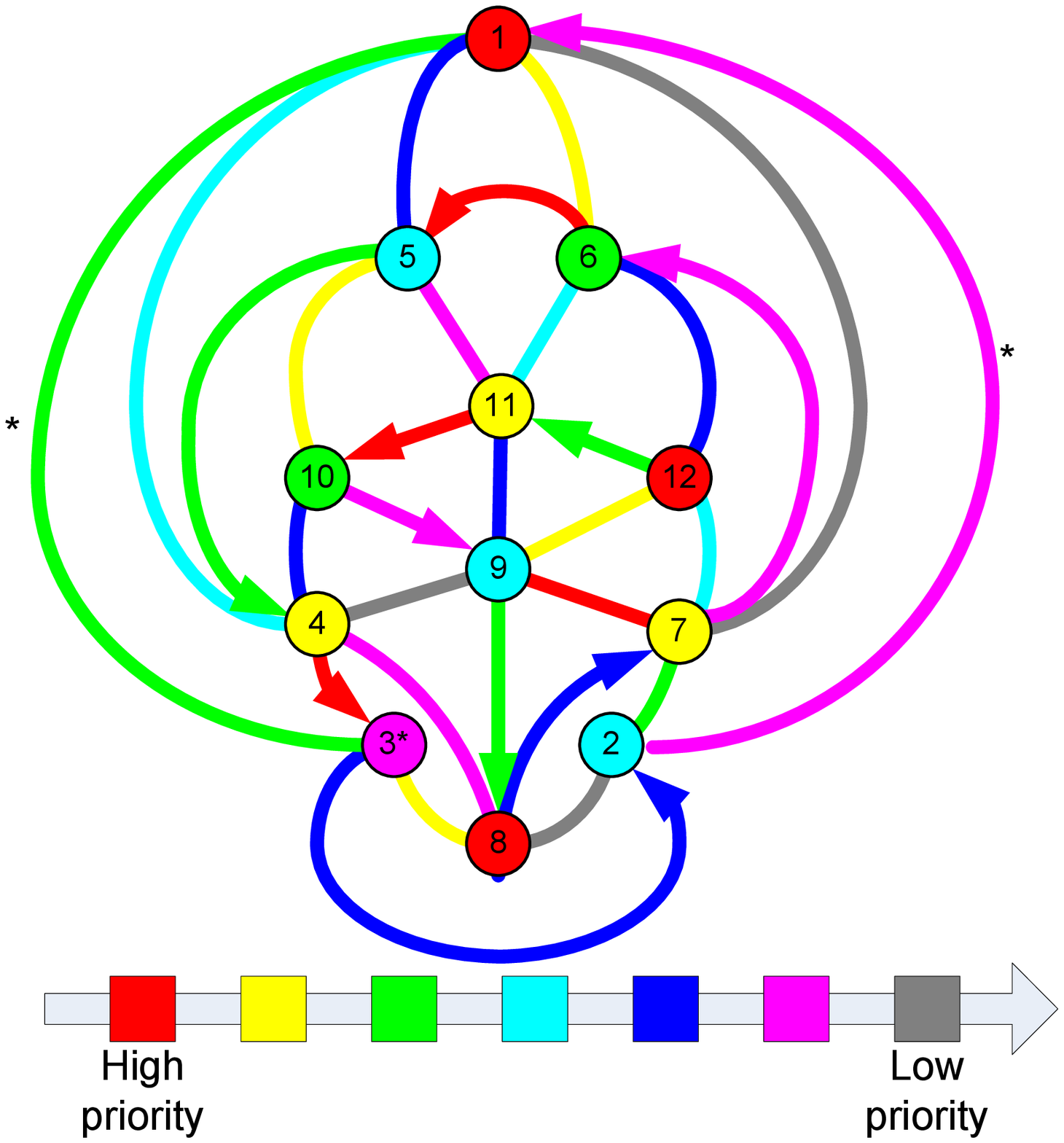}
\caption{Total coloring of an maximal planar graph with $\Delta(G)=6$.}
\label{A:fig:10}       
\end{figure}

\subsection{Spiral Chain Entire Coloring}

Kronk and Mitchem have conjectured that any plane graph of maximum degree $
\Delta $ can be colored entirely (simultaneous coloring of vertices, edges
and faces) with $\Delta (G)+4$ colors and showed that this true for $\Delta
(G)=3$ [16]. Other results on this conjecture are first
by Borodin for $\Delta (G)\geq 12$ and then $\Delta (G)\geq 7$ and finally
improved to $\Delta (G)\geq 6$ by using discharging and non-existence of an
minimal counter example by Sanders and Zhao [21]. The cases $\Delta (G)\in
\left\{ 4,5\right\} $ remain undecided. Our solution to entire coloring of
plane graph is based on the algorithms given for total and vertex spiral
chain coloring of planar graphs and valid for $\Delta (G)\geq 3.$

\textbf{Algorithm 4.}

Step 1. Find total coloring of $G$ by using spiral chain total coloring
algorithm.

Step 2. Find four coloring of the dual of $G^{\prime }$ by using vertex
spiral chain coloring algorithm. Use only the last four colors of the set $
C=\left\{ c_{1},c_{2},...,c_{(\Delta +4)}\right\}$.

The main theorem can be stated as:

\begin{theorem}
Every plane graph with maximum degree $\Delta \geq 3$ is entirely $(\Delta
+4)$-colorable.
\end{theorem}

\begin{figure}[htp]
\centering
\includegraphics[width=0.6\textwidth]{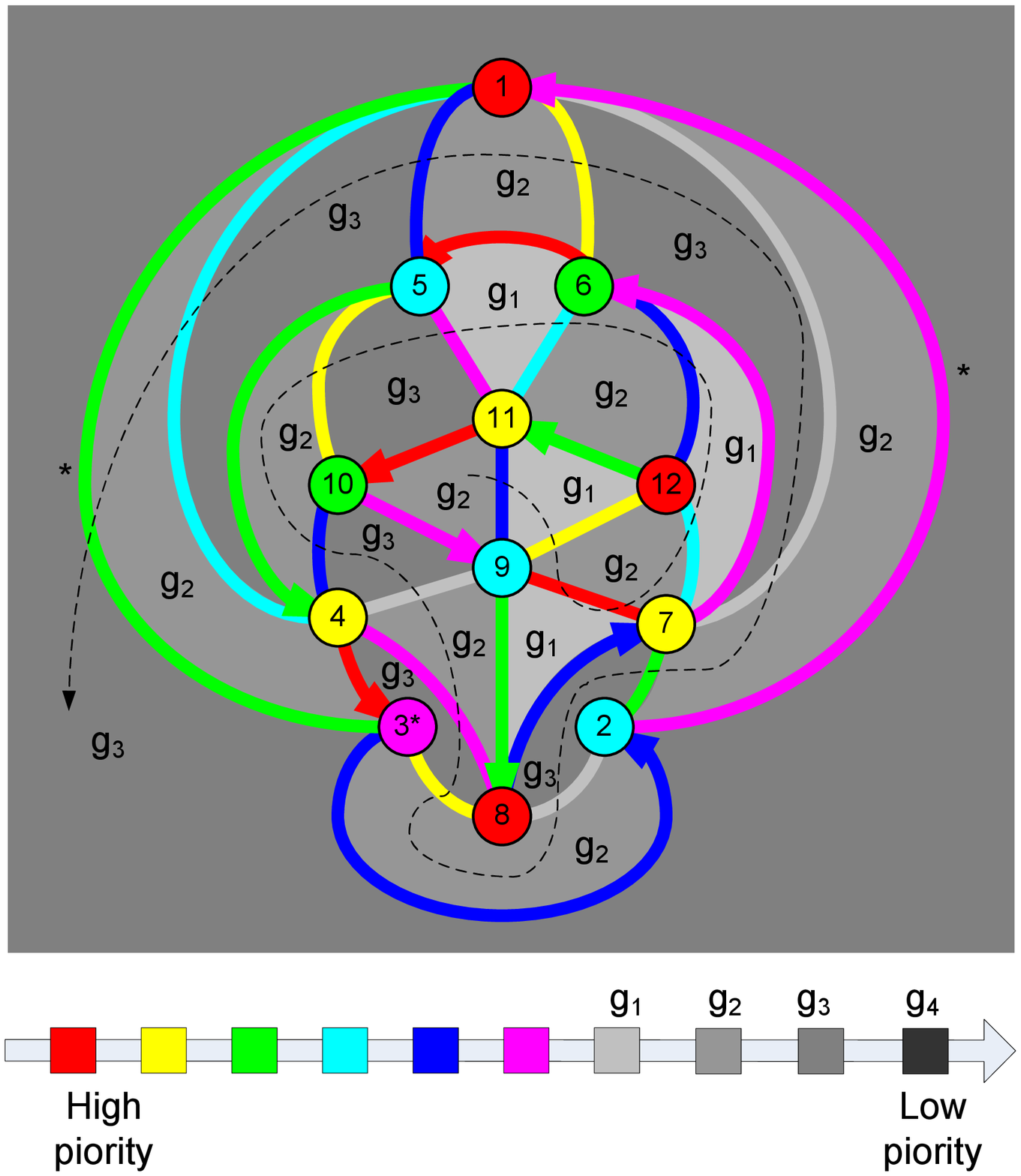}
\caption{Entire coloring of a planar graph.}
\label{A:fig:11}       
\end{figure}

Fig. 11 illustrates the spiral-chain entire coloring algorithm.

\section{Vertex Three-Colorability}

\subsection{Three color problem with triangles}

Gr\"{u}nbaum has shown that planar graphs with at most three triangles are $3$-colorable [26]. His conjecture that any planar graph having triangles apart from each other at least distance $d\geq1$ are  $3$-colorable leads to series of counterexample. Here distance $d$ between the two triangles is the length of the shortest path in the planar graph. Meinikov and Aksionov's counterexample shown in Figure 13 shows that for $d\geq3$ the graph $G$ is not $3$-colorable [25],[38]. The reason for this impasse is the $(R,Y)$-Kempe chain (shown in red dashed line) would not let us to change the only vertex colored in $R$ to a $Y$ color, adjacent to the blue $B$ colored vertex in the graph. But it is possible to reach the same conclusion that graph shown in Figure 13 can only be four colorable since it contains a unique $K_4$ as a minor; hence by the settled part of Hadwiger's conjecture it is chromatic $4$-critical. We have shown a four coloring of $G$ by using spiral chain coloring algorithm in which the last vertex colored by the fourth color blue $B$.

\begin{figure}[htp]
\centering
\includegraphics[width=0.7\textwidth]{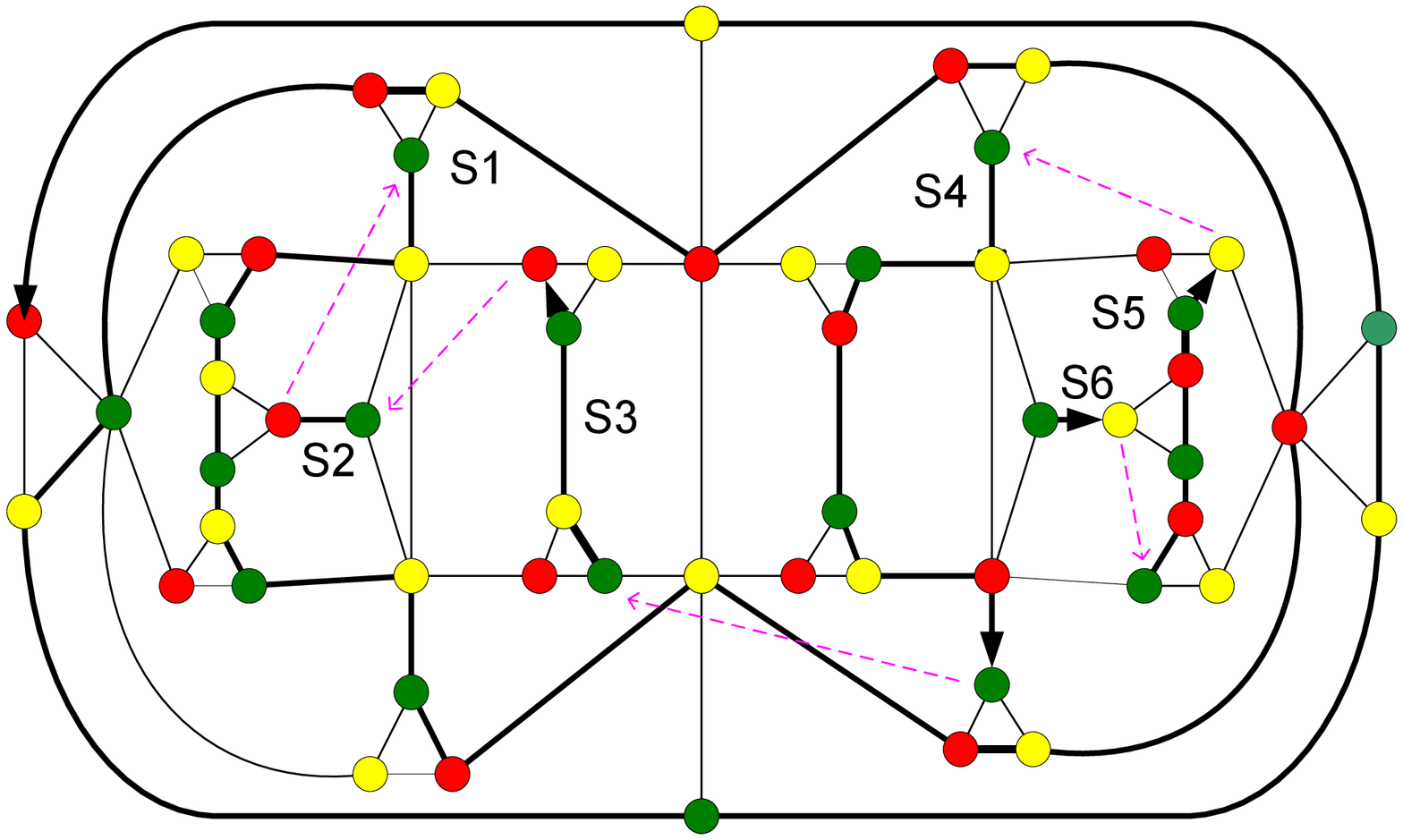}
\caption{An example for the spiral chain solution of three-colorability.}
\label{A:fig:12}       
\end{figure}

\begin{figure}[htp]
\centering
\includegraphics[width=0.8\textwidth]{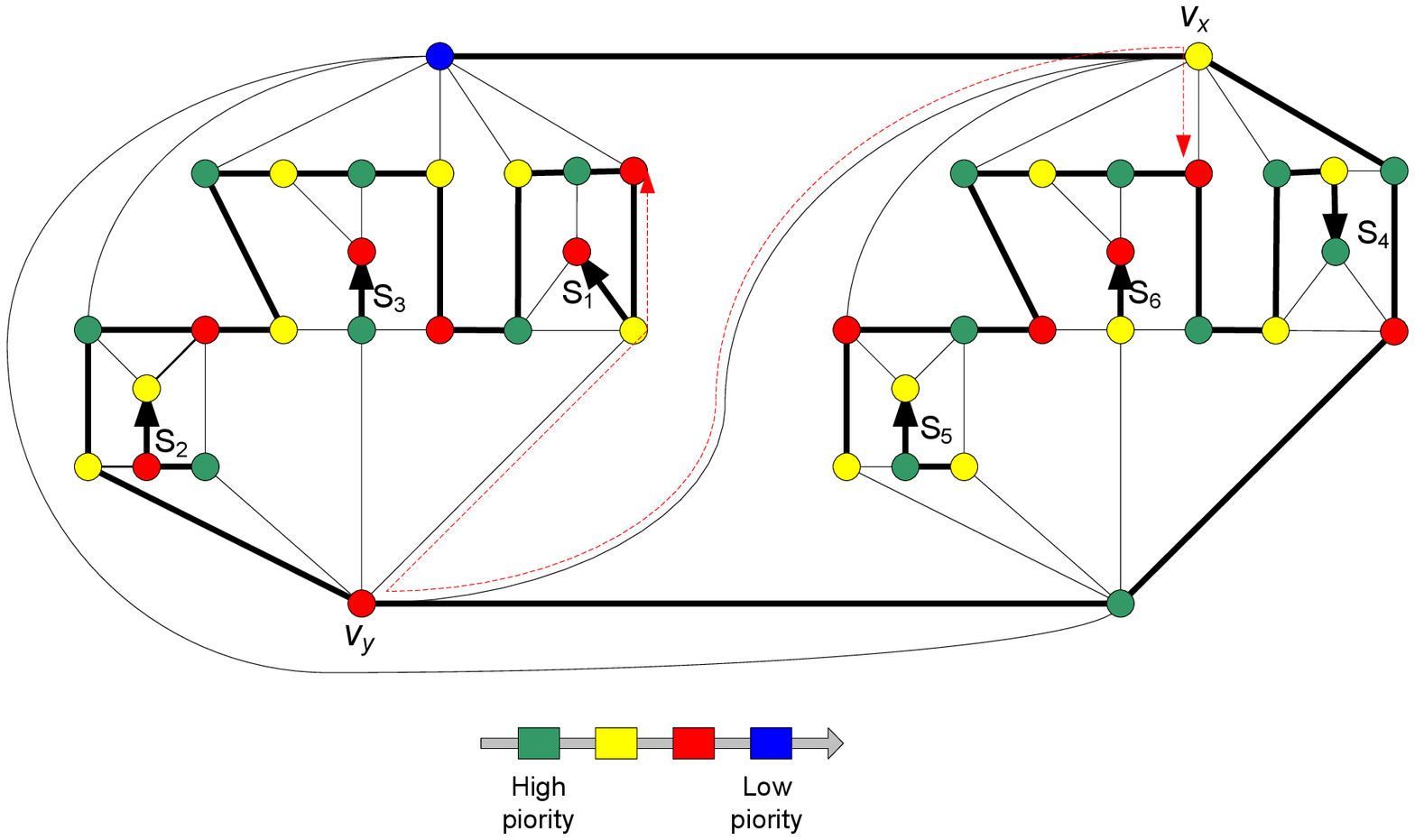}
\caption{Meinikov-Aksionov's counterexample.}
\label{A:fig:13}       
\end{figure}

In 1976 Steinberg conjectured that any planar graph without $4$-
and $5$-cycles is $3$ colorable. In [21] we have given algorithmic proof
to this conjecture based again on spiral chains. Here for the sake of
completeness we repeat the algorithm which is exactly same in principle with
the algorithms given in this paper but the size of the color set is
restricted to $3$ i.e., $C=\left\{ 1,2,3\right\} $ or $\left\{ G,Y,R\right\}
$. Let us denote by $G_{6}$, the class of planar graphs without cycles of
size from $4$ to $6$ and assume that for any vertex $v\in V(G_{6})$ we have $
\deg (v)\geq 3$. That is forbidden subgraphs in $G_{6}$ are $C_{4},C_{5}$
and any two triangles (cycle of length three) with an common edge. Assume
that we have found all spiral chains $\ S_{1},S_{2},...,S_{k}.$ Suppose that
we have completed spiral chain coloring of $G_{6}$ with only using three
colors (see [21]) and arrive at the last vertex $v_{1}$ of the spiral chain $%
S_{1}.$ Now \ any vertex on the outer-cycle $C_{o}$ of $G_{6}$ can be either
a vertex of an triangle or a non-triangle vertex. Let \ $
v_{1},v_{2},...v_{k},k\geq 6$ be the set of vertices of $C_{o}$. Let $
(v_{i}v_{i+1})\in S_{1},i=1,2,....,k-1$ and $(v_{1}v_{k})\notin S_{1}.$Let
us call a "gadget" to a subgraph consist of two triangles with an common
vertex. If $k$ is even when vertices of $C_{o}$ must be colored
alternatingly with $G$ and $Y$ starting $c(v_{1})=Y$ and if $k$ is odd all
vertices of $C_{o}$ again colored by $G$ and $Y$ except the $c(v_{k})=R$.
Let us consider a vertex $u\in C_{o}$ with $c(u)=G$. If we would join $u$
and $v_{1}$ without violating the cycle-property of $G_{6}$ and obtain a new
graph $G_{6}^{\prime }$ then since now we have new spiral chain $
S_{1}^{\prime },$ the new outer-cycle $C_{o}^{\prime }$ would be colored at
most three colors. From this we conclude that a possible counter-example to
the spiral chain coloring algorithm is the one with maximum degree at the
first vertex $v_{1}$of $S_{1}$. That is to say that is there graph $G$ with
an outer-cycle $C_{o}$ so that spiral chain coloring color vertices of $
C_{o} $ with three colors $G,Y,R$ such that the last vertex $\ v_{1}$ in the
spiral chain forcibly colored by $B$ (see for example Fig. 2 of almost three
colorable graph)? On the other hand in $G_{6}$ the outer-cycle vertices must
be in the form of serially connected gadgets or vertex of a cycle $C_{i}$
such that $\left\vert C_{i}\right\vert \geq 6$. It can be shown that in this
case $C_{o}$ can be colored with two colors and leaving room for the vertex $%
v_{1}$ to be colored with the third color.

\begin{theorem}
Planar graphs without $4$ and $5$ cycles are $3$-colorable.
\end{theorem}

\begin{figure}[htp]
\centering
\includegraphics[width=0.4\textwidth]{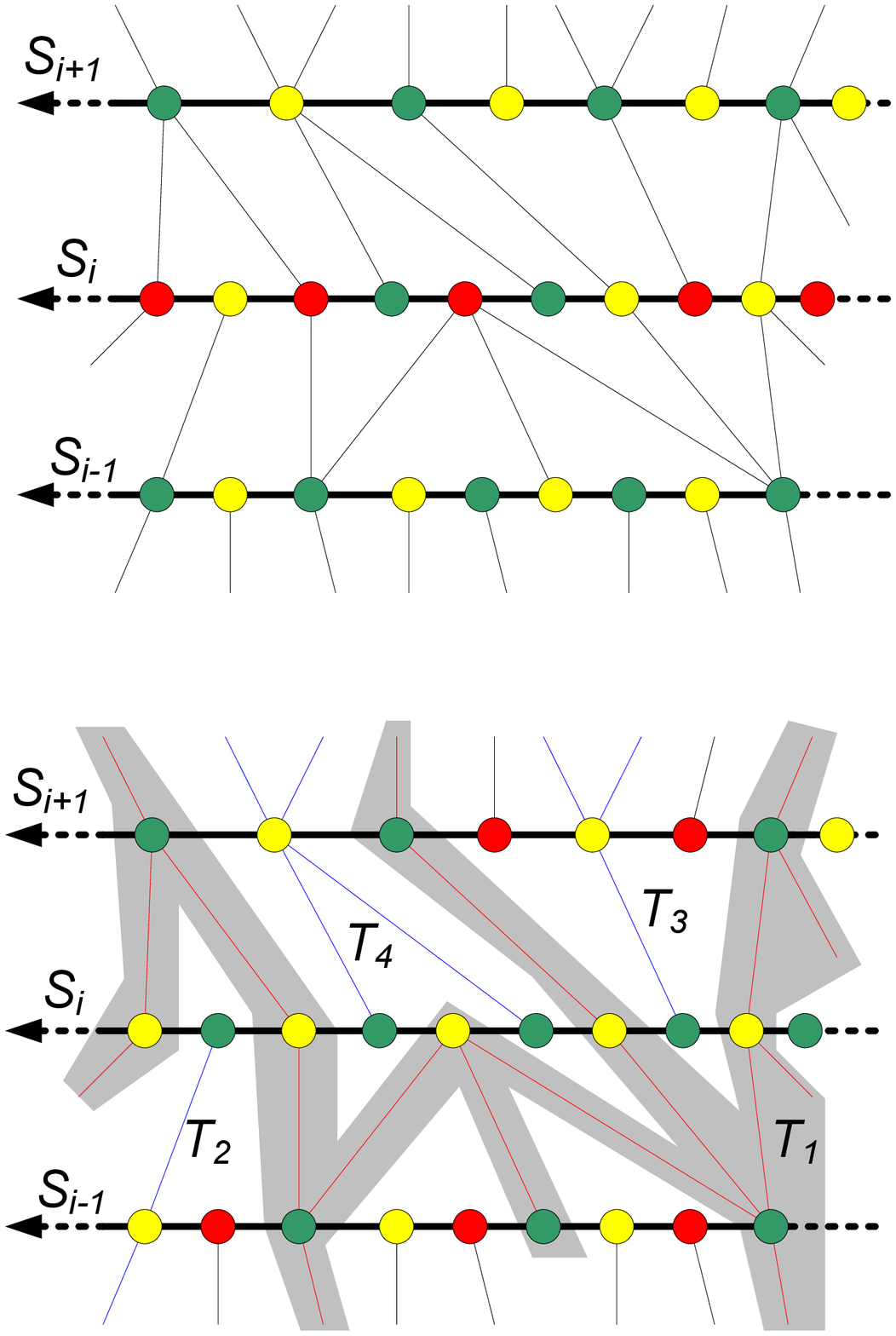}
\caption{Termination of spiral chain edge coloring algorithm with respect to the configuration at vertex $v_2$.}
\label{A:fig:6}       
\end{figure}

In [28] we investigate the $3$-colorability problem of a planar graph in
more stronger way than the above theorem under the spiral chain coloring. In
fact if $G$ has a certain type of $(\gamma ,\beta )$-sequences then it is
always possible to color it with $3$ colors by the spiral chains.

Planar graph shown in Fig.12 is taken from [27] and re-colored by the spiral
chain coloring algorithm.

\subsection{Gr\"{o}tzsch's Theorem Re-visited}
 Although planar graphs without triangles have been shown to be 3-colorable by H. Gr\"{o}etzsch [23] in 1958, on going research is still underway on this problem from the point of algorithmic complexity [35], simplification of the proof [35],[39],[40] and various counterparts on the different surfaces [36],[37]. Best bound so far is $O(nlog n)$ has been given by Kowalik [11]. In this section we will give much simpler proof to this theorem with bound O(n) by the use of spiral chain coloring algorithm [19]-[21],[28].
 It is an easy fact that chromatic number  $\chi$ of a cycle $C$ is $2$ if $|C| \equiv 0(mod2)$ and is $3$ if $|C| \equiv 1(mod2)$. Let $v_{x}$ be a vertex of a cycle $C$ and color the vertices as follow:
 (1) Start from vertex $v_{x}$ and color the vertices of $C$ in clockwise direction with the sequence of colors $1,3,1,3,...$.
 (2) Start from vertex $v_{x}$ and color the vertices of $C$ in counter clockwise direction with the sequence of colors $1,2,1,2,...$. If we start
 coloring (1) and (2) at the same time we end up with a vertex $v_{y} \in C$ for which $c(v_{x'})=2$ or $3$ if $|C| \equiv0(mod2)$. If $|C| \equiv 3(mod4)$ the two sequences end up at two adjacent vertices $v_{y'}$ and $v_{y''}$ such that $c(v_{y'})=2$ and $c(v_{y''})=3$. But color conflict arises when $|C| \equiv 1(mod4)$ such that $c(v_{y'})=c(v_{y''})=1$. This simple observation has some importance when we dealing with the algorithmic approaches e.g., spiral chain coloring to the planar graphs without triangles.

Spiral chains in $G(\ntriangleright)$
We may assume that minimum vertex degree is $\delta \geq 3$  since degree-two vertex has no effect on the three-colorability of triangle-free planar graph $G(\ntriangleright)$.
Our proof of Gr\"{o}etzsch's three color theorem is based on the following two lemmas:
\begin{lemma}
In any spiral chain decompostion of $G(\ntriangleright)$ non-spiral edges form a spanning union of cycles and trees.
\end{lemma}

\begin{lemma}
The number of third color (Red) used in the algorithm is at most equal to the number of odd cycles in $G(\ntriangleright)$.
\end{lemma}
Spiral-chain coloring algorithms
    Let $C={G,Y,R}$ be the set of three colors green, yellow and red. In the algorithm color green $(G)$ has a priority over yellow $(Y)$ and red $(R)$ and color yellow $(Y)$ has a priority over red $(R)$. Let ${S_1,S_2,...,S_{k}}$  be the set of spiral chains (assume all in clockwise directions) of $G(\ntriangleright)$. Here spiral chain (path) $S_1$ is constructed from arbitrarily selected vertex of the outer-cycle $C_{o}$ of $G( \ntriangleright)$ and continue in clockwise direction selecting all the vertices of $C_{o}$ and then continue same way to the other inner vertices (see for the details [19]-[21],[28]).

\textbf{Algorithm 5.}
Color the vertices of the set of spiral chains $S={S_{k},S_{k-1},...,S_1}$ (ordered backward direction with respect of the construction of spiral chains) using the colors of the set $C$. Coloring rule of a vertex $v\in S_{i},1\leq i\leq k$ is "use high-priority color whenever possible". If $c(v)=R$ then use $(R,G)$-Kempe chain switching or $(R,Y)$-Kempe chain switching using non-spiral-edges to re-color vertex $v$ as $c'(v)=G$ or $Y$.

\emph{Proof.} Since $G(\ntriangleright)$ is triangle-free and spiral-chains in $S$ are ordered (shelling structure) from outer spiral-chain towards an inner spiral-chain we can start coloring the innermost $S_{k}$ with two colors $G$ and  $Y$. Suppose we have arrived to coloring of a vertex $v$ with $v,u,w \in S_{k}$ such that $(vu)\in E(S_{k})$ and $(vw)\in E(S_{k})$, where $E(S_{k})$ and $E(S_{k})$ are respectively the spiral-chain and non-spiral-chain edge sets of $S_{k}$. If $c(u)=Y$ and $c(w)=G$ then use $(G,R)$-Kempe chain switching starting from vertex $w$ and go to inner colored region of $G(\ntriangleright)$. Hence we recolor $v$ as $c'(v)=G$.If $c(u)=G$ and $c(w)=Y$ then use $(Y,R)$-Kempe chain switching starting from vertex $w$ and go to inner colored region of $G(\ntriangleright)$. Hence we recolor $v$ as $c'(v)=Y$. In both cases we avoid of coloring current vertex $v$ by $R$. This means that it is possible for a full revolution of spiral-chain $S_{k}$ (called spiral-segment in [19]) all vertices can be colored by $G$ and $Y$. It is clear that we can repeat the above argument for other spiral chains. That is all red colored vertices which are unavoidable for odd cycles (by Lemma 5) are pushed into the inner spiral segments vertices of $G(\ntriangleright)$. Let $v_1$ is last vertex of $S_1$(first vertex in the construction of spiral chains). By the use of Kempe-chain switching for the outer spiral segment of $S_1$ in worst situation we color $c(v_1)=R$. Hence $G(\ntriangleright)$ is $3$-colorable.

\textbf{Algorithm 6.}
Let $S={S_{k},S_{k-1},...,S_1}$ be the set of spiral chains of $G(\ntriangleright)$. Let $F={T_{p},T_{p-1},...,T_1}$ be the forest of the set of trees formed by the non-spiral edges of $G(\ntriangleright)$. Color the vertices of trees in $F$ with green $(G)$ and yellow $(Y)$. If  $v \in T_{i},u \in T_{j},w \in T_{k}$,  $i\neq j\neq k$ with ${(vu),(vw)}\in E(S_{l})$, $1\leq  l\leq k$ such that $c(u)=G$ and $c(w)=Y$ then color $v$ as $c(v)=R$.

\emph{Proof.}Let $E(S_{i})$ and $E(S_{i})$, $i=1,2,...,k$ respectively be the sets of spiral and non-spiral edge sets of $G(\ntriangleright)$. It is not difficult to see that since $deg(v)\geq 3$ for all $v\in V(G(\ntriangleright))$, $F=\bigcup_{i=1}^{k}E(S_{i})$ is a disjoint union of trees. Denote these trees by the set $F={T_{p},T_{p-1},...,T_1}$ which is a spanning forest of $G(\ntriangleright)$. That is we have $V(T_{p})\bigcup V(T_{p-1})\bigcup...\bigcup V(T_1)=V(G(\ntriangleright)$ and $V(T_{p})\bigcap V(T_{p-1})\bigcap...\bigcap V(T_1)=\emptyset$. Since a tree is also a bipartite graph its chromatic number is two. Therefore we color vertices of each tree $T_{i}\in F$ with colors $G$ and $Y$. Then we re-color arbitrarily one of the vertex to red $R$ of all monochromatic (both green $G$ or both yellow $Y$) vertex-pairs $(uv)\in E(S_{p}),1\leq p\leq k, u\in T_{i}$ and $v\in T_{j},i\neq j$  i.e, $c(u)=c(v)=G$ or $Y\Rightarrow c'(u)=G$ and $c'(v)=R$. This results a $3$-coloring of $G(\ntriangleright)$.

\begin{theorem}
Spiral chain coloring algorithms $5$ and $6$ color any triangle-free planar graph $G(\ntriangleright)$ with three colors.
\end{theorem}

\begin{figure}[htp]
\centering
\includegraphics[width=0.6\textwidth]{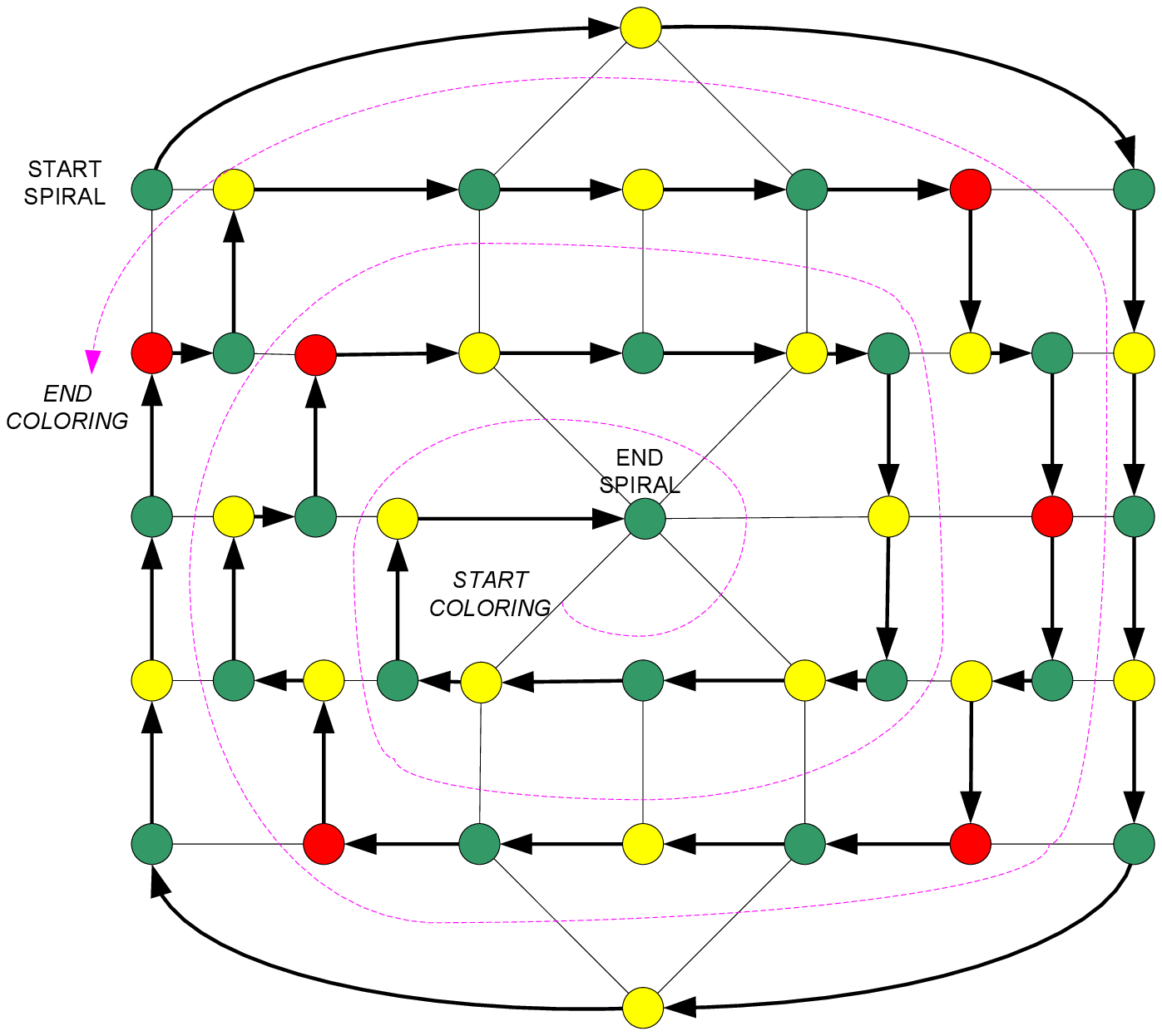}
\caption{Three coloring of a triangle-free planar graph by Algorithm 6.}
\label{A:fig:15}       
\end{figure}

Figure 1.14 illustrates the algorithms given above for three parallel spiral chains $S_{i+1}$(upper-spiral-chain), $S_{i}$ (middle-spiral-chain) and $S_{i-1}$ (lower-spiral-chain). We have also shown $3$-coloring of the trees of the non-spiral edges. In Fig. 2.15 we have illustrated spiral chain coloring algorithm for a triangle-free planar graph for which non-spiral edges induce two cycles of length $6$ and trees . Bold lines indicate spiral edges while thin dashed lines indicate non-spiral edges. In the graph there are two spiral chains $S_1$ and $S_2$.

\section{Concluding Remarks}

In this paper we have given solutions to several planar graph coloring
conjectures with the use of spiral chains. The author's 2004 algorithmic 
spiral-chain coloring proof of the famous four color theorem opens new avenues 
to the other graph coloring problems [19].

The natural question is this : 
\emph{"Why spiral-chain makes the solution of the problem so easy?"} 
Our answer to this question are several folded.
Firstly there is a famous conjecture from the complexity theory that whether
$P$ is equal $NP$ or not. It is well known result that if a problem in the
class $NP$ has been shown to be in $P$ then all other problems in the $NP$
class would have efficient solutions. Similarly we can say that many
problems related with the graphs, particularly planar graphs, would have
simple solutions if all of these graphs have Hamiltonian cycles or paths.
But we know that some of graphs are not Hamiltonian and in fact finding one
in a giving graph is not easy. Hamilton path problem is difficult but the
algorithmic answer of finding spiral-chain in graph is almost trivially very
easy. Moreover we can easily decided when and how the graph has more than
one spiral chains. Therefore spiral chain would act as navigator and paths
decomposition in the graph coloring for us to reach the solution. In other
words spiral chain is a road-map for efficient coloring algorithm.

Secondly the use of the spiral-chain reduces the number of the cases
considerably in the proof. Many other proof methods in the graph coloring
theorems are to show nonexistence of minimal counter-examples. But this in
most of the times is a very complex task and sometimes we need to
investigate case-by-case by only using a computer. Just consider how the
possible impasse in the spiral chain proof of the four color theorem is
ruled out by re-coloring certain vertex pair in the "sailing boat" subgraph
of the maximal planar graph.

I think the third one is the most important. Suppose we start to color the
vertices of a planar graph by using spiral chain then you cannot say
beforehand whether this process partitioned the graph into, say two parts at
the end. Lastly when we color a vertex in the spiral chain we are sure that
we will consider another vertex later on that adjacent to the previously
colored vertex. That is the main idea that prevents us to fall into the
troubles like one of the most elegant "proof" in mathematics [14].

The next question about the use of spiral chains is the following:

\emph{How one can apply spiral chains to the coloring problems for non-planar graphs?}
I think a well-defined spiral chains decomposition of a complete graph may help to
devise algorithmic solution to the Hadwiger's conjecture [15] which asserts that that every 
loopless graph not contractible to the complete graph on $t+1$ vertices is $t$-colorable. When $t = 3$ this is easy, and
when $t = 4$, Wagner's theorem of 1937 shows the conjecture to be equivalent to the
four-color conjecture (the 4CC)[42]. The case $t = 5$ it is also equivalent to the 4CC. 
Without assuming the 4CC Robertson, Seymour and Thomas have shown that every minimal counterexample to Hadwiger's
conjecture when $t = 5$ is \emph{apex}, that is, it consists of a planar graph with one
additional vertex. Consequently, the 4CC implies Hadwiger's conjecture when $t = 5$,
because it implies that apex graphs are $5$-colorable [43].

\end{document}